\documentclass[letterpaper,11pt]{article}

\usepackage{ucs}
\usepackage[utf8x]{inputenc}
\usepackage{graphicx}
\usepackage{amsfonts}
\usepackage{dsfont}
\usepackage{amssymb}
\usepackage{amsmath}
\usepackage{amsthm}
\usepackage{enumerate}
\usepackage{stmaryrd}
\usepackage{fullpage}
\usepackage{ifthen}
\usepackage{subfigure}
\usepackage{epic}
\usepackage{authblk}
\usepackage{textcomp}

\usepackage[hypertexnames=false,colorlinks=true,linkcolor=blue,citecolor=blue]{hyperref}
\usepackage[numbers,comma,square,sort&compress]{natbib}
\usepackage[letterpaper,text={7in,9in},centering]{geometry}

\usepackage{color}
\usepackage{titlesec}
\graphicspath{{eps/}{pdf/}}

\newcommand{\bqq}{\begin{equation}}
\newcommand{\eqq}{\end{equation}}
\newcommand{\bqs}{\begin{equation*}}
\newcommand{\eqs}{\end{equation*}}

\newcommand{\R}{\mathbb{R}} 
\newcommand{\C}{\mathbb{C}}

\newcommand{\rmO}{\mathrm{O}}

\newcommand{\A}{\mathcal{A}}

\newcommand{\cC}{\mathcal{C}}

\newcommand{\F}{\mathcal{F}}
\newcommand{\G}{\mathcal{G}}
\newcommand{\cH}{\mathcal{H}}

\newcommand{\J}{\mathcal{J}}
\newcommand{\K}{\mathcal{K}}
\newcommand{\cl}{\mathcal{L}}
\newcommand{\M}{\mathcal{M}}
\newcommand{\cN}{\mathcal{N}}
\newcommand{\cP}{\mathcal{P}}
\newcommand{\T}{\mathcal{T}}
\newcommand{\W}{\mathcal{W}}
\newcommand{\X}{\mathcal{X}}
\newcommand{\Y}{\mathcal{Y}}

\newcommand{\mG}{\mathbf{G}}
\newcommand{\mH}{\mathbf{H}}

\numberwithin{equation}{section}

\newtheorem{lem}{Lemma}[section]
\newtheorem{thm}{Theorem}
\newtheorem{prop}[lem]{Proposition}
\newtheorem{cor}[lem]{Corollary}
\newtheorem{rmk}[lem]{Remark}

\newtheorem{hyp}[lem]{Hypothesis}

\makeatletter\@addtoreset{figure}{section}\makeatother

\makeatletter \@addtoreset{equation}{section} \makeatother

\newcommand{\etal}{\textit{et al.}\ }

\newenvironment{Proof}%
 {\begin{trivlist} \item[]{\bf Proof. }}%
{\hspace*{\fill}$\rule{.4\baselineskip}{.4\baselineskip}$\end{trivlist}}

\title{Fredholm properties of nonlocal differential operators via spectral flow}
\author[1]{Gr\'egory Faye}
\author[2]{Arnd Scheel}
\affil[1,2]{\small University of Minnesota,
School of Mathematics,
206 Church Street S.E.,
Minneapolis, MN 55455, USA}

\begin{document}
\maketitle

\begin{abstract}
We establish Fredholm properties for a class of nonlocal differential operators. Using mild convergence and localization conditions on the nonlocal terms, we also show how to compute Fredholm indices via a generalized spectral flow, using crossing numbers of generalized spatial eigenvalues. We illustrate possible applications of the results in a nonlinear and a linear setting. We first prove the existence of small viscous shock waves in nonlocal conservation laws with small spatially localized source terms. We also show how our results can be used to study edge bifurcations in eigenvalue problems using Lyapunov-Schmidt reduction instead of a Gap Lemma.
\end{abstract}

{\noindent \bf Keywords:} Nonlocal operator; Fredholm index; Spectral flow; Nonlocal conservation law; Edge bifurcations.\\


\section{Introduction}

\subsection{Motivation}

Our aim in this paper is the study of the following class of nonlocal linear operators:
\bqq
\label{mapTintro}
\T: H^1(\R,\R^n) \longrightarrow L^2(\R,\R^n),\quad  U\longmapsto \frac{d}{d\xi}U-\widetilde{\K}_\xi \ast U
\eqq
where the matrix convolution kernel $\widetilde{\K}_\xi(\zeta)=\widetilde{\K}(\zeta;\xi)$ acts via
\bqs
\widetilde{\K}_\xi \ast U(\xi)=\int_\R \widetilde{\K}(\xi-\xi';\xi)U(\xi')d\xi'.
\eqs
Operators such as \eqref{mapTintro} appear when linearizing at coherent structures such as traveling fronts or pulses in nonlinear nonlocal differential equations. One is interested in properties of the linearization when analyzing robustness, stability or interactions of these coherent structures. A prototypical example are neural field equations which are used in mathematical neuroscience to model cortical traveling waves. They typically take the form \cite{pinto-ermentrout:01}
\begin{subequations}
\label{eq:nfePop}
\begin{align}
\partial_t u(x,t)&=-u(x,t)+\int_{\R}K(|x-x'|)S(u(x',t))dx'-\gamma v(x,t) \\
\partial_t v(x,t)&=\epsilon(u(x,t)-v(x,t))
\end{align}
\end{subequations}
for $x\in\R$ and with $\gamma$, $\epsilon$ positive parameters.  The nonlinearity $S$ is the firing rate function and the kernel $K$ is often referred to as the connectivity function. It encodes how neurons located at position $x$ interact with neurons located at position $x'$ across the cortex. The first equation describes the evolution of the synaptic current $u(x,t)$ in the presence of linear adaptation which takes the form of a recovery variable $v(x,t)$ evolving according to the second equation. In the moving frame $\xi=x-ct$, equations \eqref{eq:nfePop} can be written as
\begin{subequations}
\label{eq:nfePopMF}
\begin{align}
\partial_t u(\xi,t)&=c \partial_\xi u(\xi,t)-u(\xi,t)+\int_{\R}K(|\xi-\xi'|)S(u(\xi',t))d\xi'-\gamma v(\xi,t) \\
\partial_t v(\xi,t)&=c \partial_\xi v(\xi,t)+\epsilon(u(\xi,t)-v(\xi,t)),
\end{align}
\end{subequations}
such that stationary solutions $(u(\xi),v(\xi))$ satisfy
\begin{subequations}
\label{eq:nfePopTW}
\begin{align}
-c \frac{d}{d\xi}u(\xi)&=-u(\xi)+\int_{\R}K(|\xi-\xi'|)S(u(\xi'))d\xi'-\gamma v(\xi) \\
-c \frac{d}{d\xi}v(\xi)&=\epsilon(u(\xi)-v(\xi)).
\end{align}
\end{subequations}
The linearization of \eqref{eq:nfePopMF} at a particular solution $(u_0(\xi),v_0(\xi))$ of \eqref{eq:nfePopTW} takes the form
\begin{subequations}
\label{eq:nfePopTWlin}
\begin{align}
\partial_t u(\xi,t)&=c \partial_\xi u(\xi,t)-u(\xi,t)+\int_{\R}K(|\xi-\xi'|)S'(u_0(\xi'))u(\xi',t)d\xi'-\gamma v(\xi,t) \\
\partial_t v(\xi,t)&=c \partial_\xi v(\xi,t)+\epsilon(u(\xi,t)-v(\xi,t)).
\end{align}
\end{subequations}
Denoting $U=(u,v)$ and $\mathcal{L}_0$ the right-hand side of \eqref{eq:nfePopTWlin}, the eigenvalue problem associated with the linearization of \eqref{eq:nfePopMF} at $(u_0,v_0)$ reads 
\bqq
\label{eigvalpbm}
\lambda U = \mathcal{L}_0 U.
\eqq
This eigenvalue problem can be cast as a first-order nonlocal differential equation
\bqq
\label{eigvalpbmcast}
\frac{d}{d\xi}U(\xi)=\widetilde{\K}_\xi^\lambda\ast U(\xi)
\eqq
where 
\bqs
\widetilde{\K}_\xi^\lambda(\zeta)=-\frac{1}{c}\left(
\begin{matrix}
-(1+\lambda)\delta_0+K(|\zeta|)S'(u_*(\xi-\zeta)) & -\gamma \delta_0\\
\epsilon \delta_0 & -(\epsilon+\lambda)\delta_0
\end{matrix}
\right)
\eqs
and $\delta_0$ denotes the Dirac delta at $0$.

The differential systems \eqref{eq:nfePopTW} and  \eqref{eigvalpbmcast} can be viewed as  systems of functional differential equations of mixed type since the convolutional term  introduces both advanced and retarded terms. Such equations are notoriously difficult to analyze. Our goal here is threefold. First, we establish Fredholm properties of such operators. Second we give algorithms for computing Fredholm indices. Last, we show how such Fredholm properties can be used to analyze perturbation and stability problems.

For local differential equations, a variety of techniques is available to study such problems. For example, in the case of the Fitzhugh-Nagumo equations, written in moving frame $\xi=x-ct$,
\begin{subequations}
\label{eq:FHN}
\begin{align}
\partial_t u&=c \partial_\xi u+\partial_{\xi\xi} u +f(u)-\gamma v \\
\partial_t v&=c \partial_\xi v + \epsilon(u-v)
\end{align}
\end{subequations}
with a bistable nonlinearity $f$, spectral properties of the linearization of \eqref{eq:FHN} at a stationary solution $(u_*(\xi),v_*(\xi))$
\bqs
\cl_* := \left(\begin{matrix} c \partial_\xi +\partial_{\xi\xi}+f'(u_*) & -\gamma \\ \epsilon & c\partial_\xi -\epsilon \end{matrix} \right),
\eqs
are encoded in exponential dichotomies of the first-order equation \cite{palmer:88,sandstede:02}
\bqq
\label{eq:1order}
\frac{d}{d\xi}U(\xi)=A(\xi,\lambda)U(\xi), \quad A(\xi,\lambda)=\left(\begin{matrix} 0 & 1 & 0 \\ \lambda - f'(u_*) & -c & -\gamma \\ -\frac{\epsilon}{c} & 0 &  \frac{\lambda+\epsilon}{c} \end{matrix} \right).
\eqq
In particular, $\cl^*-\lambda$ is a Fredholm operator if and only if \eqref{eq:1order} has exponential dichotomies on $\R^-$ and $\R^+$. Unfortunately, for nonlocal equations \eqref{eigvalpbmcast}, neither existence of exponential dichotomies nor Fredholm properties are known in general. Spectral properties of nonlocal operators such as $\T$ in  \eqref{mapTintro} are understood mostly in the cases where $\T-\lambda$ is Fredholm with index zero and $U$ is scalar. We mention the early work of Ermentrout \& McLeod \cite{ermentrout-mcleod:93} who proved that the Fredholm index at a traveling front is zero in the case where $\gamma=0$ (no adaptation) for the neural field system \eqref{eq:nfePop}. Using comparison principles, De Masi \etal proved stability results for traveling fronts in nonlocal equations arising in Ising systems with Glauber dynamics and Kac potentials \cite{demasi-etal:95}. In a more general setting, yet relying on comparison principles, Chen \cite{chen:97} showed the existence and asymptotic 
stability of traveling fronts for a class of nonlocal equations, including the models studied by Ermentrout \& McLeod and De Masi \etal. Bates \etal \cite {bates-etal:97}, using monotonicity and a homotopy argument, also studied the existence, uniqueness, and  stability of traveling wave solutions in a bistable, nonlinear, nonlocal equation. 

More general results are available when the interaction kernel is a finite sum of Dirac delta measures. In particular, the interaction kernel has finite range in that case. Such interaction kernels arise in the study of lattice dynamical systems. Mallet-Paret established Fredholm properties and showed how to compute the Fredholm index via a spectral flow \cite{mallet-paret:99}. His methods are reminiscent of Robbin \& Salamon's work \cite{robbinsalamon:95}, who established similar results for operators $\frac{d}{d\xi}+A(\xi)$ where $A(\xi)$ is self-adjoint but does not necessarily generate a semi-group. For the operators studied in \cite{mallet-paret:99}, Fredholm properties are in fact equivalent to the existence of exponential dichotomies for an appropriate formulation of \eqref{mapTintro} as an infinite- dimensional evolution problem \cite{mallet-paret-verduyn-lunel:01,harterich-etal:02}.

Our approach extends Mallet-Paret's results \cite{mallet-paret:99} to infinite-range kernels. We do not know if a dynamical systems formulation in the spirit of \cite{mallet-paret-verduyn-lunel:01,harterich-etal:02} is possible. Our methods blend some of the tools in \cite{robbinsalamon:95} with techniques from \cite{mallet-paret:99}. In the remainder of the introduction, we give a precise statement of assumptions and our main results.

\subsection{Main results --- summary}

We are interested in proving  Fredholm properties for 
\bqs
\T:U\longmapsto \frac{d}{d\xi}U-\widetilde{\K}_\xi \ast U.
\eqs
Our main results assume the following properties for $\widetilde{\K}_\xi$
\begin{itemize}
\item {\bf Exponential localization}: the kernel $\widetilde{\K}_\xi$ is exponentially localized, uniformly in $\xi$; see Section \ref{subsecnotations}, Hypotheses \ref{hypK} and \ref{hypA}.
\item {\bf Asymptotically constant}: there exist constant kernels $\widetilde{\K}^\pm$ such that $\widetilde{\K}_\xi\underset{\xi\rightarrow\pm \infty}{\longrightarrow}\widetilde{\K}^\pm$; see Section \ref{subsecnotations}, Hypotheses \ref{hypK} and \ref{hypA}.
\item {\bf Asymptotic hyperbolicity}: the asymptotic kernels $\widetilde{\K}^\pm$ are hyperbolic; see Section \ref{subsecnotations}, Hypothesis \ref{zero}, and Section \ref{subsecasympt}.
\item {\bf Asymptotic regularity}: the complex extensions of the Fourier transforms of $\widetilde{\K}^\pm$ are bounded and analytic in a strip containing the imaginary axis; see Section \ref{subsecnotations}, Hypothesis \ref{bound}.
\end{itemize} 

Our main results can then be summarized as follows.

\begin{thm}
Assume that the interaction kernel $\widetilde{\K}_\xi$ satisfies the following properties: exponential localization, asymptotically constant, asymptotic hyperbolicity, and asymptotic regularity. Then the nonlocal operator $\T$ defined in \eqref{mapTintro} is Fredholm on $L^2(\R)$ and its index can be computed via its spectral flow.
\label{thmFormal}
\end{thm}

As a first example, we study shocks in nonlocal conservation laws with small localized sources of the form
\bqq
\label{nvclintro}
U_t =\left( \K \ast F(U) +  G(U) \right)_x + \epsilon H(x,U,U_x),\quad U\in \R^n.
\eqq
Similar types of conservation laws have been studied in \cite{chmaj:07,du-etal:12}. More precisely, using a monotone iteration scheme, Chmaj proved the existence of traveling wave solutions for \eqref{nvclintro} with $\epsilon=0$, $U\in\R$, \cite{chmaj:07}. Du \etal proposed to study nonlocal conservation laws more systemically and described interesting behavior in the inviscid nonlocal Burgers' equation \cite{du-etal:12}. We show how our results can help study properties of shocks in such systems \eqref{nvclintro}. We prove that for small localized external sources there exist small undercompressive shocks of index $-1$, that is, $\# \left\{ \text{outgoing characteristics}\right\}$ $=$ $\# \left\{ \text{ingoing characteristics}\right\} $. Shocks can be parametrized by values on ingoing "characteristics" in the case when characteristic speeds do not vanish. For vanishing characteristic speeds, we show the existence of undercompressive shocks with index $-2$, that is, $\# \left\{ \text{outgoing characteristics}\right\}$ $=$ $\# 
\left\{ \text{ingoing characteristics}\right\} +2$. Here, we use the term characteristic informally, a precise definition via the dispersion relation is given in Section \ref{secNCL}.

As a second example, we consider bifurcation of eigenvalues from the edge of the essential spectrum. It has been recognized early \cite{simon:76} that localized perturbations of operators can cause eigenvalues to emerge from the essential spectrum. More recently, spatial dynamics methods have helped to treat a much larger class of eigenvalue problems using analytic extensions of the Evans function into the essential spectrum, thus tracking eigenvalues into and beyond the essential spectrum; see \cite{gz,ks}. This extension, usually referred to as the Gap Lemma, was used to track stability and instability in a conservation law during spatial homotopies \cite{pogan-scheel:10,pogan-scheel:11}, without referring to spatial dynamics but rather to a local tracking function constructed via Lyapunov-Schmidt and matching proceedures. In Section \ref{secGL}, we will show that such an approach is possible for nonlocal equations, using the Fredholm properties established in our main results.

\subsection{Set-up of the problem}

We are interested in studying linear nonlocal differential equations that can be written as:
\bqq
\label{eqNL}
\frac{d}{d\xi}U(\xi)=\int_\R \K(\xi-\xi';\xi)U(\xi')d\xi'+\sum_{j\in \J} A_j(\xi) U(\xi-\xi_j)+H(\xi).
\eqq
Here $U(\xi),H(\xi)\in\C^n$, and $\K(\zeta;\xi),A_j(\xi)\in\M_n(\C)$, $n\geq 1$, the space of $n\times n$ complex matrices. The set $\J$ is countable and the shits $\xi_j$ satisfy (without loss of generality)
\bqq
\label{shifts}
\xi_1=0,\quad \xi_j\neq\xi_k, \quad j \neq k \in \J.
\eqq

For each $\xi \in \R$, we define $\A(\xi)$ by
\bqq
\label{opA}
\A(\xi):=\left( \K(~ \cdot ~ ;\xi), \left(A_j(\xi) \right)_{j\in\J}\right),
\eqq
such that we may write \eqref{eqNL} as
\bqq
\label{eqNLm}
\frac{d}{d\xi}U(\xi)= \cN[\A(\xi)]\cdot U(\xi) + H(\xi),
\eqq
where $\cN[\A(\xi)]$ denotes the linear nonlocal operator
\bqq
\label{opNA}
\cN[\A(\xi)]\cdot U(\xi):=\int_\R \K(\xi-\xi';\xi)U(\xi')d\xi'+\sum_{j\in \J} A_j(\xi) U(\xi-\xi_j).
\eqq
We denote $\K_\xi:=\K(~\cdot~;\xi)$ and write \eqref{opNA}  as a generalized convolution \bqq
\label{opNAr}
\cN[\A(\xi)]\cdot U =\left[ \K_\xi  + \sum_{j\in\J} A_j(\xi) \delta_{\xi_j}\right] \ast U .
\eqq
Here $\ast$ refers to convolution on $\R$
\bqs
(W_1 \ast W_2)(\xi)=\int_\R W_1(\xi-\xi')W_2(\xi')d\xi',
\eqs
and $\delta_{\xi_j}$ is the Dirac delta at $\xi_j \in \R$. 

Setting $H\equiv 0$, we obtain the homogeneous system
\bqq
\label{eqNLhom}
\frac{d}{d\xi}U(\xi)= \cN[\A(\xi)]\cdot U(\xi).
\eqq

A special case of \eqref{opNAr} are constant coefficient operators $\A(\xi)$
\bqs
\A(\xi)=\left(\K^0(~\cdot~), \left(A_j^0 \right)_{j\in\J}\right):=\A^0, \quad \forall\,\xi\in\R.
\eqs
We have
\bqq
\label{opNAconstant}
\cN[\A^0]\cdot U =\left[ \K^0  + \sum_{j\in\J} A_j^0 \delta_{\xi_j}\right] \ast U 
\eqq
and
\bqq
\label{eqNLconstant}
U'(\xi)= \cN[\A^0]\cdot U(\xi).
\eqq

Associated with \eqref{eqNLhom}, we have the linear operator 
\bqq
\label{opT}
\mathcal{T_\A}:=\dfrac{d}{d\xi}-\cN[\A(\xi)].
\eqq 

\subsection{Notations and hypotheses}\label{subsecnotations}

We denote by $\cH$ and $\W$ the Hilbert spaces $L^2(\R,\C^n)$ and $H^1(\R,\C^n)$ equipped with their usual norm
\bqs
\| U \|_\cH := \underset{k=1\cdots n}{\max}\| U_k \|_{L^2(\R)},
\eqs
and
\bqs
\| U \| _\W := \| U' \|_\cH+\| U \|_\cH.
\eqs

For a function $\K_\xi=\K(~\cdot~;\xi):\R\rightarrow L^1_\eta(\R,\M_n (\C))$, $\eta>0$, we define its norm as
\bqs
|| \K_\xi ||_\eta:= \underset{(k,l)\in  \llbracket 1,n \rrbracket^2}{\max}\|\K_{k,l}(~\cdot~;\xi)e^{\eta |~\cdot ~ |}\|_{L^1(\R)}.
\eqs
We also introduce the following norm for the kernel $\K\in\cC^1\left(\R,L^1_\eta(\R,\M_n (\C))\right)$, 
\bqs
||| \K |||_{\infty,\eta} := \underset{\xi\in\R}{\sup}  \left\| \K_\xi \right\|_\eta + \underset{\xi\in\R}{\sup}  \left\| \frac{d}{d\xi}\K_\xi \right\|_\eta.
\eqs
For a function $A\in \cC^1(\R,\M_n(\C))$ we define its norm as
\bqs
\| A \|_n:= \underset{\xi\in\R}{\sup}  \left\|A(\xi) \right\|_{\M_n(\C)} + \underset{\xi\in\R}{\sup}  \left\| \frac{d}{d\xi}A(\xi) \right\|_{\M_n(\C)}.
\eqs

Finally we denote by $\tau$ the linear transformation that acts on $\K_\xi$ as $\tau \cdot \K_\xi:=\K(~\cdot~;~\cdot~+\xi) $ and we naturally define $\tau\cdot\K:\xi \longmapsto \tau \cdot \K_\xi$. We can now give further assumptions on the maps $\K$ and $\left(A_j\right)_{j\in\J}$.

\begin{hyp}\label{hypK}
There exists $\eta>0$ such that the matrix kernel $\K$ satisfies the following properties:
\begin{enumerate}
\item  $\K$ belongs to $\cC^1\left(\R,L^1_\eta(\R,\M_n (\C))\right)$;
\item $\K$ is localized, that is,
\begin{subequations}
\label{etaK}
\begin{align}
||| \K |||_{\infty,\eta}&<\infty~,\\
 ||| \tau \cdot \K |||_{\infty,\eta}&<\infty~;
\end{align}
\end{subequations}
\item there exist two functions $\K^\pm\in L^1(\R, \M_n (\C))$ such that
\bqq
\label{limK}
\underset{\xi \rightarrow \pm \infty}{\lim}\K(\zeta;\xi)=\K^\pm(\zeta)
\eqq
uniformly in $\zeta\in\R$ and
\begin{subequations}
\label{limJeta}
\begin{align}
\underset{\xi \rightarrow \pm \infty}{\lim} \| \K_\xi - \K^\pm \|_\eta&=0\\
\underset{\xi \rightarrow \pm \infty}{\lim} \| \tau\cdot \K_\xi - \K^\pm \|_\eta&=0.
\end{align}
\end{subequations}
\end{enumerate}
\end{hyp}

\begin{hyp}\label{hypA}
The matrices $A_j$ satisfy the properties:
\begin{enumerate}
\item $A_j \in \cC^1\left(\R,\M _n(\C)\right)$ for all $j\in\J$;
\item with $\eta$ defined in Hypothesis \ref{hypK}, we have,
\bqq
\label{etaA}
\sum_{j\in\J}\|A_j\|_ne^{\eta|\xi_j|}<\infty~;
\eqq
\item there exist $A_j^\pm\in \M_n (\C)$ such that
\bqq
\label{limA}
\underset{\xi \rightarrow \pm \infty}{\lim}A_j(\xi)=A_j^\pm, \quad \sum_{j\in\J} \|A_j^\pm\|_{\M_n (\C)}e^{\eta|\xi_j|}  <\infty, \quad  j\in\J
\eqq
and
\bqq
\label{finitesum}
\underset{\xi \rightarrow \pm \infty}{\lim} \sum_{j\in\J} \|A_j(\xi)-A_j^\pm\|_{\M_n (\C)} e^{\eta|\xi_j|}=0. 
\eqq
\end{enumerate}
\end{hyp}

Note that if we define the map $\A$ as
\bqq
\label{mapa}
\begin{matrix}
\A: & \R & \longrightarrow &  L^1_\eta (\R,\M_n(\C)) \times \ell^1_\eta \left(\M_n(\C) \right)  \\
& \xi & \longmapsto & \A(\xi)=\left( \K(~ \cdot ~ ;\xi), \left(A_j(\xi) \right)_{j\in\J}\right)
\end{matrix}
\eqq
then, when Hypotheses \ref{hypK} and \ref{hypA} are satisfied, $\A\in\cC^1(\R,L^1_\eta (\R,\M_n(\C)) \times \ell^1_\eta \left(\M_n(\C) \right))$ and is bounded. Here we have implicitly defined
\bqs
\ell^1_\eta \left(\M_n(\C) \right)=\left\{ (A_j)_{j\in\J} \in \M_n(\C)^\J ~|~ \sum_{j\in\J}\|A_j\|_{\M_n (\C)}e^{\eta|\xi_j|}<\infty \right\}.
\eqs

\begin{hyp}\label{zero}
We assume that for all $\ell \in \R$
\bqq
\label{det}
d^\pm(i\ell):=\det \left( i\ell~\mathbb{I}_n -\widehat{\K^\pm}(i \ell)-\sum_{j\in\J}A_j^\pm  e^{-i\ell \xi_j} \right) \neq 0
\eqq
where $\widehat{\K^\pm}$ are the complex Fourier transforms of $\K^\pm$ defined by
\bqs
\widehat{\K^\pm}(i\ell) = \int_\R \K^\pm(\xi)e^{-i\ell \xi}d\xi.
\eqs
\end{hyp}

\begin{hyp}\label{bound}
We assume that, with the same $\eta>0$ as in Hypotheses \ref{hypK} and \ref{hypA}, the complex Fourier transforms 
\bqs
\nu\longmapsto \widehat{\K^\pm}(\nu)+\sum_{j\in\J}A_j^\pm  e^{-\nu \xi_j}
\eqs 
extend to bounded analytic functions in the strip $\mathcal{S}_\eta:=\left\{\nu \in \C~|~ |\Re(\nu)| < \eta \right\}$.
\end{hyp}

\subsection{Main results}

We can now restate our informal Theorem \ref{thmFormal} which we split in two separate theorems. The first theorem states the Fredholm property of the nonlocal operator $\T_\A$ while the second gives a characterization of the Fredholm index via the spectral flow.

\begin{thm}[The Fredholm Alternative]
\label{thmFredholm}
Suppose that Hypotheses \ref{hypK}, \ref{hypA}, and \ref{zero} are satisfied. Then the operator $\T_\A:\W\rightarrow\cH$ is Fredholm. Furthermore, the Fredholm index of $\T_\A$ depends only on the limiting operators $\A^\pm$, the limits of $\A(\xi)$ as $\xi\rightarrow\pm\infty$. We denote $\iota(\A^-,\A^+)$ the Fredholm index $\mathrm{ind}~\T_\A$.

\end{thm}

\begin{thm}[Spectral Flow Theorem]
\label{thmFlow}
Assume that Hypotheses \ref{hypK}, \ref{hypA}, \ref{zero}, and \ref{bound} are satisfied and suppose, further, that there are only finitely many values of $\xi_0\in\R$ for which $\A(\xi_0)$ is not hyperbolic. Then the Fredholm index of $\T_\A$
\bqq
\label{index}
\iota(\A^-,\A^+)=-\mathrm{cross}(\A)
\eqq 
is the net number of roots of \eqref{eqNLhom} which cross the imaginary axis from left to right as $\xi$ is increased from $-\infty$ to $+\infty$; see Section \ref{subsecCrossings} for a precise definition.
\end{thm}

\begin{rmk}
Similar Fredholm results hold for higher-order differential operators with nonlocal terms. This can be seen by transforming into a system of first-order equations, or, more directly, by following the proof below, which treats the main part of the equation as a generalized operator pencil, thus allowing for more general forms of the equation.
\end{rmk}

\paragraph{Outline.} This paper is organized as follows. We start in Section \ref{PrelimNot} by introducing some notation and basic material needed in the subsequent sections. Section \ref{fredholm} is devoted to the proof of Theorem \ref{thmFredholm} while in Section \ref{flow} we prove Theorem \ref{thmFlow}. Finally in Section \ref{application}, we apply our results to nonlocal conservation laws with spatially localized source term and to nonlocal eigenvalue problems with small spatially localized nonlocal perturbations.

\section{Preliminaries and notation}\label{PrelimNot}

Consider Banach spaces  $\X$ and $\Y$. We let $\cl(\X,\Y)$ denote the Banach space of bounded linear operators $\T:\X\rightarrow\Y$, and we denote the operator norm by $\|\T\|_{\cl(\X,\Y)}$. We write $\text{rg }\T$ for the range of $\T$ and $\text{ker }\T$ for its kernel,
\bqs
\text{rg }\T:=\left\{ \T U \in \Y ~;~U\in \X \right\}\subset  \Y,
\qquad
\text{ker }\T:=\left\{U\in\X ~;~ \T U =0 \right\}\subset\X.
\eqs
In the proof of Theorem \ref{thmFredholm}, we shall use the following Lemma; see \cite{schwarz:93} for a proof.

\begin{lem}[Abstract Closed Range Lemma]\label{ACRL}
Suppose that $\cal X$, $\cal Y$ and $\cal Z$ are Banach spaces, that $\cal T : \cal X \rightarrow \cal Y$ is a bounded linear operator, and that $\cal R : \cal X \rightarrow \cal Z$ is a compact linear operator. Assume that there exists a constant $c>0$ such that
\bqs
\| U \|_{\cal X} \leq c \left(\| \T U\|_{\cal Y}+\|\mathcal{R} U\|_{\cal Z} \right), \quad \forall\, U \in \cal X.
\eqs
Then $\cal T$ has closed range and finite-dimensional kernel. 
\end{lem}

Let us recall that a bounded operator $\T:\X\rightarrow\Y$ is a Fredholm operator if
\begin{itemize}
\item[(i)] its kernel $\text{ker }\T$ is finite-dimensional;
\item[(ii)] its range $\text{rg }\T$ is closed; and
\item[(iii)] $\text{rg }\T$ has finite codimension.
\end{itemize}
For such an operator, the integer
\bqs
\text{ind~}\T:=\dim \left(\text{ker }\T\right) -\text{codim}\left(\text{rg }\T\right)
\eqs
is called the Fredholm index of $\T$.

\subsection{Adjoint equation}

We introduce the formal adjoint equation of \eqref{eqNLhom} as
\bqq
\label{adjoint}
\frac{d}{d\xi}U(\xi):=\cN [\A(\xi)]^*\cdot U(\xi)=-\int_\R \K^*(\xi'-\xi;\xi')U(\xi')d\xi'-\sum_{j\in \J} A_j^*(\xi+\xi_j) U(\xi+\xi_j)
\eqq
with $\K^*$ and $A_j^*$ denoting the conjugate transposes of the matrices $\K$ and $A_j$, respectively. Elementary calculations give  that $\cN [\A(\xi)]^*=\cN[\widetilde{\A}(\xi)]$ where
\bqs
\widetilde{\A}(\xi)=\left(\widetilde{K}(~\cdot~;\xi),(\widetilde{A}_j(\xi))_{j\in\J} \right)
\eqs
and $\widetilde{K}$ and $\widetilde{A}_j$ are defined as
\begin{align*}
\widetilde{K}(\zeta;\xi)&=-\K^*(-\zeta;-\zeta+\xi)\quad \forall\, \zeta\in\R, \\
\widetilde{A}_j(\xi)&=-A_j^*(\xi+\xi_j)\quad \forall\, j\in\J.
\end{align*}
Note that $\widetilde{K}$ and $\widetilde{A}_j$ also satisfy Hypotheses \ref{hypK} and  \ref{hypA}.

Considering $\T_\A$ as a closed, densely defined operator on $\cH$, we find that the adjoint $\T_\A^*:\W\subset\cH\rightarrow\cH$ is given through
\bqq
\label{adjT}
\T_\A^*=-\dfrac{d}{d\xi}+\cN [\A(\xi)]^*.
\eqq

\subsection{Asymptotically autonomous systems}\label{subsecasympt}

Associated to the constant coefficient system \eqref{eqNLconstant}  is the characteristic equation 
\bqq
\label{eqcharacteristic}
d^0(\nu):=\det \Delta_{\A^0}(\nu)=0
\eqq
where
\bqq
\label{matDconstant}
\Delta_{\A^0}(\nu)= \nu~\mathbb{I}_n -\widehat{\K^0}(\nu)-\sum_{j\in\J}A_j^0  e^{-\nu \xi_j}, \quad \nu\in\C.
\eqq
Note that the characteristic equation possesses imaginary roots precisely when there exist solutions of the form $e^{i\ell\xi}$ to \eqref{eqNLconstant}. More generally, roots of $d^-(\nu)$ detect pure exponential solutions to \eqref{eqNLconstant}. 
We say that this constant coefficient system is hyperbolic when
\bqq
\label{hyperbolicity}
d^0(i\ell)\neq 0,\quad  \forall\, \ell \in \R.
\eqq
In the specific case considered here, when $\widehat{\K^0}$ is a bounded analytic function in the strip $\mathcal{S}_\eta$, there are only finitely many roots of \eqref{eqcharacteristic} in the strip. One can think of roots $\nu$ of \eqref{eqcharacteristic} as generalized eigenvalues to the generalized eigenvalue problem \eqref{opNAconstant}.

We say that the system \eqref{eqNLhom} is asymptotically autonomous at $\xi=+\infty$ if 
\bqs
\underset{\xi\rightarrow +\infty}{\lim}\A(\xi)=\A^+
\eqs
where $\A^+$ is constant. In this case, of course, \eqref{eqNLconstant} with $\A^0=\A^+$ is called the limiting equation at $+\infty$. If in addition, the limiting equation is hyperbolic, then we say that \eqref{eqNLhom} asymptotically hyperbolic at $+\infty$. We analogously define asymptotically autonomous and asymptotically hyperbolic at $-\infty$. If \eqref{eqNLhom} is asymptotically autonomous at both $\pm\infty$, we simply say that \eqref{eqNLhom} is asymptotically autonomous, asymptotically hyperbolic if  asymptotically hyperbolic at $\pm\infty$.

In the case of the constant coefficient system \eqref{eqNLconstant} it is straightforward to see that we have 
\bqs
\Delta_{\A^{0*}}(\nu)=-\Delta_{\A^0}(-\bar\nu)^*,
\eqs
so that
\bqs
\det \Delta_{\A^{0*}}(\nu)=(-1)^n \det \Delta_{\A^0}(-\nu).
\eqs
This implies that system \eqref{eqNLconstant} is hyperbolic if and only if its adjoint is hyperbolic.

\section{Fredholm properties}\label{fredholm}

For each $T>0$, we define $\cH(T)=L^2([-T,T],\C^n)$ and $\W(T)=H^1([-T,T],\C^n)$. It is easy to see that the inclusion $\W(T)\hookrightarrow \cH(T)$ defines a compact operator such that the restriction operator 
\begin{align*}
\mathcal{R}&:\W\rightarrow \cH(T)\\
&~~ U \mapsto U_{[-T,T]}
\end{align*}
 is a compact linear operator and $\|\mathcal{R} U\|_{\cH(T)}=\|U\|_{\cH(T)}$.
\begin{lem}\label{l:rs}
There exist constants $c>0$ and $T>0$ such that 
\bqq
\label{estimateT}
\| U \|_\W \leq c \left( \|U\|_{\cH(T)}+\| \T_\A U\|_{\cH}\right)
\eqq
for every $U\in \W$.
\end{lem}
\begin{Proof}
Following \cite{robbinsalamon:95}, we divide the proof into three steps.

\underline{\textit{Step - 1}}\quad  For each $U\in\W$, we have
\bqs
\| \T_\A U \|_\cH = \left\| \frac{d}{d\xi}U(\xi)-\cN[\mathcal{A}(\xi)]\cdot U \right\|_\cH \geq \left\| \frac{d}{d\xi}U \right\|_\cH- C\| U \|_\cH,
\eqs
where the constant $C>0$ can be chosen as
\bqs
C= n \left( \sqrt{|||\K |||_{\infty,\eta}|||\tau\cdot \K |||_{\infty,\eta}} +\sum_{j\in\J} \| A_j \|_n\right).
\eqs
Indeed, fix $k\in\llbracket1,n\rrbracket$, and estimate
\begin{align*}
\int_\R \left| \left(\K_\xi\ast U\right)_k(\xi) \right|^2d\xi &\leq n \sum_{l=1}^n\int_\R\left(\int_\R \left|  \K_{k,l}(\xi-\xi';\xi)U_l(\xi') \right|d\xi'\right)^2 d\xi \\
&\leq n \sum_{l=1}^n\int_\R\left(\int_\R \left|  \K_{k,l}(\xi-\xi';\xi)\right|^{1/2} \left|  \K_{k,l}(\xi-\xi';\xi)\right|^{1/2} \left| U_l(\xi') \right|d\xi'\right)^2 d\xi \\
&\leq n \sum_{l=1}^n\int_\R \left( \int_\R \left|  \K_{k,l}(\xi-\xi';\xi)\right| d\xi' \right)  \left( \int_\R \left|  \K_{k,l}(\xi-\xi';\xi)\right|  \left| U_l(\xi') \right|^2d\xi' \right)  d\xi \\
&\leq n ||| \K |||_{\infty,\eta}  \sum_{l=1}^n \int_\R \left(\int_\R \left|  \K_{k,l}(\xi-\xi';\xi)\right|d\xi \right) \left| U_l(\xi') \right|^2 d\xi' \\
& \leq n^2 ||| \K |||_{\infty,\eta} ||| \tau\cdot \K |||_{\infty,\eta} \| U \|_{\cH}^2.
\end{align*}
Similarly, one obtains
\bqs
\int_\R \left| \left(A_j(\xi) U(\xi-\xi_j)\right)_k \right|^2 d\xi \leq n^2 \| A_j \|_n^2 \| U \|_{\cH}^2.
\eqs

This proves the estimate \eqref{estimateT} with $T=\infty$:
\bqq
\label{estTinf}
\| U \|_\W \leq c_1 \left( \|U\|_{\cH}+\| \T_\A U\|_{\cH}\right).
\eqq

\underline{\textit{Step - 2}}\quad In the second step, we prove the estimate for a hyperbolic, constant coefficient system \eqref{eqNLconstant}, 
\bqs
\cN[\A^0]\cdot U=  \left[ \K^0  + \sum_{j\in\J} A_j^0 \delta_{\xi_j}\right] \ast U.
\eqs
Applying Fourier transform to $f=\T_{\A^0} U$ gives
\bqs
\left(i\ell \mathbb{I}_n  - \widehat{\K^0}(i\ell)- \sum_{j\in\J}A_j^0 e^{-i\ell\xi_j} \right)\widehat{U}(i \ell)=\widehat{f}(i\ell) \quad \forall\, \ell\in\R.
\eqs
Using the fact that $\A^0$ is hyperbolic ($d^0(i\ell) \neq 0$), we can invert
\bqs
\widehat{U}(i \ell)=\left(i\ell \mathbb{I}_n  - \widehat{\K^0}(i\ell)- \sum_{j\in\J}A_j^0 e^{-i\ell\xi_j} \right)^{-1} \widehat{f}(i\ell) \quad \forall\, \ell\in\R.
\eqs
This implies that
\bqs
\| \widehat{U}\|_{\cH} \leq \underset{\ell \in \R}{\sup}\left(i\ell \mathbb{I}_n  - \widehat{\K^0}(i\ell)- \sum_{j\in\J}A_j^0 e^{-i\ell\xi_j} \right)^{-1} \|\widehat{f}\|_\cH,
\eqs
and, using the Fourier-Plancherel theorem, we obtain
\bqs
\| U \|_{\cH} \leq c_1  \| \T_{\A^0}U  \|_\cH \quad \forall\, U \in \W,
\eqs
for some constant $c_1>0$. Using the first step, we finally have the inequality
\bqq
\label{eq_inter}
\| U \|_\W \leq c_2  \| \T_{\A^0}U  \|_\cH \quad \forall\, U \in \W,
\eqq
with $c_2>0$.

\underline{\textit{Step - 3}}\quad We want to prove that there exist $T>0$ such that, if $U(\xi)=0$ for $|\xi| \leq T-1$, $U\in\W$, we have
\bqq
\label{ineqstep3}
\| U\|_W \leq c_3  \|\T_\A U\|_\cH.
\eqq
To do so, we first prove that inequality \eqref{ineqstep3} is satisfied for functions $U^\pm\in \W$, of the form
\bqq
\label{Upm}
U^+(\xi)= 0 \text{ for } \xi \leq T-1 \text{ and } U^-(\xi)=0 \text{ for } \xi \geq -T+1.
\eqq
We remark that Hypotheses \ref{hypK} and \ref{hypA} ensure the existence of $T>0$ and $\epsilon(T)>0$ such that, if $U^\pm\in\W$ are defined as above, the following estimates are satisfied 
\begin{subequations}\label{epsilon}
\begin{align}
\left\| \left( \K^\pm-\K_\xi \right) \ast U^\pm \right\|_\cH &\leq \frac{\epsilon(T)}{2} \|U^\pm \|_\cH,\\
\left\| \sum_{j\in\J} \left( A_j^\pm-A_j(\xi)\right)\left( \delta_{\xi_j}\ast U^\pm\right) \right\|_\cH &\leq \frac{\epsilon(T)}{2} \|U^\pm\|_\cH.
\end{align}
\end{subequations}

This ensures that for every $U^\pm\in\W$ satisfying \eqref{Upm}, we have
\bqs
\frac{1}{c_2}\| U^\pm  \|_\W \leq \| \T_{\A^\pm} U^\pm \|_\cH \leq \epsilon(T) \|U^\pm\|_\cH + \|\T_\A U^\pm\|_\cH,
\eqs
which proves inequality \eqref{ineqstep3} in that case. Here, we have used the implicit notations
\begin{align*}
\T_{\A^\pm}&=\frac{d}{d\xi}-\cN[\A^\pm],\\
\cN[\A^\pm]\cdot U^\pm&= \left[ \K^\pm  + \sum_{j\in\J} A_j^\pm \delta_{\xi_j}\right] \ast U .
\end{align*}

Finally, if $U\in\W$ is such that $U(\xi)=0$ for $|\xi| \leq T-1$, we decompose $U$ as the sum $U^++U^-$, setting
\bqs
U^+(\xi)=\left\{ 
\begin{matrix}
U(\xi), & \xi \geq 0 \\
0, & \xi < 0
\end{matrix}
\right., \qquad \qquad 
U^-(\xi)=\left\{ 
\begin{matrix}
0, & \xi > 0\\
U(\xi), & \xi \leq 0
\end{matrix}
\right..
\eqs
Of course, $U^\pm$ now satisfy \eqref{Upm} and we have
\bqs
\|U\|^2_\W=\|U^+\|^2_\W+\|U^-\|^2_\W \leq c_3^2 \left( \|\T_\A U^+\|_\cH^2+\|\T_\A U^-\|_\cH^2  \right)=c^2_3\|\T_\A U\|_\cH^2,
\eqs
which gives the desired inequality.

\underline{\textit{Step - 4}}\quad Finally, the estimate \eqref{estimateT} is proved by a patching argument. We choose a smooth cutoff function $\chi:\R\rightarrow [0,1]$ such that $\chi(\xi)=0$ for $|\xi| \geq T$ and $\chi(\xi)=1$ for $|\xi| \leq T-1$. Using estimate \eqref{estTinf} for $\chi U$ and \eqref{ineqstep3} for $(1-\chi)U$, we have
\begin{align*}
\|U\|_\W & \leq \|\chi U\|_\W+\|(1-\chi)U\|_\W\\
&\leq c_1(\|\chi U\|_\cH+\|\T_\A(\chi U)\|_\cH)+c_3\|\T_\A[(1-\chi)U]\|_\cH\\
&\leq c \left((\|U\|_{\cH(T)}+\|\T_\A(U)\|_\cH\right).
\end{align*}
\end{Proof}

Together with the abstract closed range Lemma \ref{ACRL}, Lemma \ref{l:rs} immediately implies the semi-Fredholm properties for $\T_\A$ and its adjoint.
\begin{cor} 
Both, $\T_\A$ and $\T_\A^*$, considered as operators from $\W$ into $\cH$, possess closed range and finite-dimensional kernel.
\end{cor}

\begin{Proof}
We only need to verify that the Hypotheses \ref{hypK}, \ref{hypA} and \ref{zero} are satisfied for the adjoint operator $\T_\A^*$. We recall that in that case we have
\bqs
\T_\A^*=-\frac{d}{d\xi}+\cN[\widetilde{\A}(\xi)]
\eqs
where 
\bqs
\widetilde{\A}(\xi)=\left(\widetilde{K}(~\cdot~;\xi),(\widetilde{A}_j(\xi))_{j\in\J} \right)
\eqs
and $\widetilde{K}$ and $\widetilde{A}_j$ are defined as
\begin{align*}
\widetilde{K}(\zeta;\xi)&=-\K^*(-\zeta;-\zeta+\xi)\quad \forall\, \zeta\in\R, \\
\widetilde{A}_j(\xi)&=-A_j^*(\xi+\xi_j)\quad \forall\, j\in\J.
\end{align*}
As a consequence, Hypotheses \ref{hypK} and \ref{hypA} are satisfied for the adjoint. Hypothesis \ref{zero} refers to asymptotic hyperbolicity of $\T_\A$. We already noticed that $\A^\pm$ is hyperbolic if and only if its adjoint $\A^{\pm*}$ is hyperbolic, which implies that Hypothesis \ref{zero} is also satisfied for the adjoint equation. By Lemma \ref{l:rs}, $\T_\A^*$ then has closed range and finite-dimensional kernel.
\end{Proof}
\begin{Proof}[of Theorem \ref{thmFredholm}] 
The above corollary implies that $\T_\A:\W \rightarrow \cH$ has finite-dimensional kernel, closed range, and finite-dimensional co-kernel given by the kernel of its adjoint $\T_\A^*$. 

To prove that the Fredholm index depends only on the limiting operators $\A^\pm$ we consider two families of operators $\A_0(\xi)$ and $\A_1(\xi)$ that satisfy  Hypotheses \ref{hypK}, \ref{hypA} and \ref{zero} with coefficients
\bqs
\A_0(\xi)=\left(\K_0(~\cdot~;\xi),  \left(A_{j,0}(\xi) \right)_{j\in\J}\right),\quad \A_1(\xi)=\left(\K_1(~\cdot~;\xi),  \left(A_{j,1}(\xi) \right)_{j\in\J}\right)
\eqs
and the same shifts $\xi_j$. We assume that the limiting operators at $\pm\infty$ are equal, that is,
\bqs
\A^\pm_0=\A^\pm_1,
\eqs
where
\bqs
\A^\pm_\sigma= \left( \K^\pm_\sigma, \left(A_{j,\sigma}^\pm\right)_{j\in\J}\right)=\underset{\xi \rightarrow \pm \xi}{\lim} \A_\sigma(\xi), \quad \sigma=0,1.
\eqs
For $0\leq \sigma \leq 1$, we define $\A_\sigma(\xi)=(1-\sigma)\A_0(\xi)+\sigma\A_1(\xi)$. Then for each such $\sigma$, $\A_\sigma$ satisfies Hypotheses \ref{hypK}, \ref{hypA} and \ref{zero} and $\T_{\A_\sigma}$ is a Fredholm operator and $\T_{\A_\sigma}$ varies continuously in $\cl(\W,\cH)$ with $\sigma$. Thus the Fredholm index of $\T_{\A_\sigma}$ is independent of $\sigma$ and only depends on the limiting operators $\A^\pm$.

\end{Proof}

\begin{rmk}
The proof immediately generalizes to a set-up where $\cH$ and $\W$ are $L^p$-based, with $1<p<\infty$, with the exception of invertibility of the asymptotic, constant-coefficient operators, where we used Fourier transform as an isomorphism. On the other hand, analyticity of the Fourier multiplier shows that the inverse is in fact represented by a convolution with an exponentially localized kernel, which gives a bounded inverse in $L^p$, so that our theorem holds in $L^p$-based spaces as well.
\end{rmk}

\begin{cor}[Cocycle property]\label{cor:cocycle}
Suppose that $\A^0,\A^1$ and $\A^2$ are hyperbolic constant coefficient operators in $L^1_\eta(\R,\M_n(\C))\times\ell^1_\eta\left(\M_n(\C) \right)$, then we have
\bqs
\iota(\A^0,\A^1)+\iota(\A^1,\A^2)=\iota(\A^0,\A^2).
\eqs
\end{cor}
\begin{Proof}
We consider, for $0\leq\sigma\leq1$, the system
\bqs
\mathcal{U}'(\xi)=\cN[\A_\sigma(\xi)]\mathcal{U}(\xi),\quad \mathcal{U}(\xi)\in\C^{2n}
\eqs
where $\A_\sigma(\xi)=\left(\K_\sigma(~\cdot~;\xi),  \left(A_{j,\sigma}(\xi) \right)_{j\in\J}\right)\in L^1_\eta(\R,\M_{2n}(\C))\times\ell^1_\eta\left(\M_{2n}(\C) \right)$
\begin{align*}
\K_\sigma(~\cdot~;\xi)&=\chi_-(\xi) \left(\begin{matrix} \K^0(~\cdot~) & 0 \\ 0 & \K^1(~\cdot~) \end{matrix}\right)+\chi_+(\xi) R(\sigma) \left( \begin{matrix} \K^1(~\cdot~) & 0 \\ 0 & \K^2(~\cdot~) \end{matrix}\right)  R(-\sigma),\\
A_{j,\sigma}\xi)&=\chi_-(\xi) \left(\begin{matrix} A_j^0 & 0 \\ 0 & A_j^1\end{matrix}\right)+\chi_+(\xi) R(\sigma) \left( \begin{matrix} A_j^1 & 0 \\ 0 & A_j^2 \end{matrix}\right)  R(-\sigma),\\
R(\sigma)&=\left(\begin{matrix} \cos\left(\frac{\pi\sigma}{2}\right) & - \sin\left(\frac{\pi\sigma}{2}\right) \\ \sin\left(\frac{\pi\sigma}{2}\right) & \cos\left(\frac{\pi\sigma}{2}\right)\end{matrix} \right)
\end{align*}
with $\chi_\pm(\xi)=(1+\tanh(\pm\xi))/2$. For all $0\leq\sigma\leq1$, $\A_\sigma(\xi)$ is asymptotically hyperbolic and satisfies Hypotheses \ref{hypK} and \ref{hypA}, thus $\T_{\A_\sigma}$ is Fredholm and the Fredholm index of $\T_{\A_\sigma}$ is independent of $\sigma$. Namely, we have $\text{ind}~ \T_{\A_{\sigma=0}}= \text{ind}~\T_{\A_{\sigma=1}}$. At $\sigma=0$ and $\sigma=1$, the equation $\mathcal{U}'(\xi)=\cN[\A_\sigma(\xi)]\mathcal{U}(\xi)$ decouples and one finds that
\begin{align*}
\text{ind}~ \T_{\A_{\sigma=0}}&=\iota(\A^0,\A^1)+\iota(\A^1,\A^2),\\
\text{ind}~\T_{\A_{\sigma=1}}&=\iota(\A^0,\A^2)+\iota(\A^1,\A^1)=\iota(\A^0,\A^2).
\end{align*}
This concludes the proof.
\end{Proof}

\section{Spectral flow}\label{flow}

Throughout this section we fix the shifts $\xi_j$. For $\rho\in\R$, we denote by $\A^\rho$ a continuously varying one-parameter family of constant coefficient operators of the form: 
\bqq
\label{eq:1pf}
\begin{matrix}
\A: & \R & \longrightarrow &  L^1_\eta(\R,\M_n(\C))\times\ell^1_\eta\left(\M_n(\C) \right) \\
& \rho & \longmapsto & \A^\rho=\left(\K^\rho(~\cdot~),  \left(A_j^\rho \right)_{j\in\J}\right).
\end{matrix}
\eqq
For simplicity, we identify the family $\A^\rho$ with its associated constant nonlocal operator $\cN[\A^\rho]$. In this section we will prove the following result which automatically gives the result of Theorem \ref{thmFlow}.

\begin{thm}\label{thmSFI}
Let $\A^\rho$, for $\rho \in\R$, a continuously varying one-parameter family of constant coefficient operators of the form \eqref{eq:1pf}. We suppose that:
\begin{itemize}
\item[(i)] the limit operators $\A^\pm$ are hyperbolic in the sense that $\forall\, \ell \in \R$
\bqs
d^{\pm}(i\ell)=\det\left(i\ell \mathbb{I}_n - \widehat{\K^{\pm}}(i\ell)-\sum_{j\in\J}A_j^{\pm} e^{-i\ell \xi_j}\right)\neq 0 ,
\eqs
\item[(ii)] $\Delta_{\A^{\rho}}(\nu)$ defined in \eqref{matDconstant} is a bounded analytic function in the strip $\mathcal{S}_\eta=\left\{\lambda \in \C~|~ |\Re(\lambda)| < \eta \right\}$ for each $\rho\in\R$. 
\item[(iii)] there are finitely many values of $\rho$ for which $\A^\rho$ is not hyperbolic.
\end{itemize}
Then
\bqq
\label{index}
\iota(\A^{-},\A^{+})=-\mathrm{cross}(\A)
\eqq 
is the net number of roots of \eqref{eqNLhom} which cross the imaginary axis from left to right as $\rho$ is increased from $-\infty$ to $+\infty$.
\end{thm}

In our approach to the proof , we approximate the family $\A^\rho$ of Theorem \ref{thmSFI} with a generic family \cite{mallet-paret:99,robbinsalamon:95}. To do so, we need to introduce some notations. We denote by $\cP:=\cP(\R,L^1_\eta (\R,\M_n(\C)) \times \ell^1_\eta\left(\M_n(\C) \right))$ the Banach space of all continuous paths for which conditions $(i)$ and $(ii)$ of Theorem \ref{thmSFI} are satisfied. And finally, define the open set 
$\cP^1:=\cC^1(\R,L^1_\eta (\R,\M_n(\C)) \times \ell^1_\eta\left(\M_n(\C) \right)\cap \cP.$

\subsection{Crossings}\label{subsecCrossings}

For any continuous path $\A$ of the form \eqref{eq:1pf}, a \textbf{crossing} for $\A$ is a real number $\rho_0$ for which $\A^{\rho_0}$ is not hyperbolic and we let
\bqs
\text{NH}(\A):=\left\{ \rho\in \R ~|~ \text{ equation } \eqref{eqNLhom} \text{ with constant coefficients } \A^\rho \text{ is not hyperbolic} \right\},
\eqs
be the set of all crossings of $\A$. Thus $\A$ satisfies condition $(iii)$ of Theorem \ref{thmSFI} if and only if $\A$ has finitely many crossings. In that case, $\text{NH}(\A)$ is a finite set that we denote by $\text{NH}(\A)=\left\{\rho_1,\dots,\rho_m \right\}$. Note that for all $\A\in\cP$ and at any crossing $\rho_0$,  the equation 
\bqs
d_{\rho_0}(\nu):=\det(\Delta_{\A^{\rho_0}}(\nu))= 0 
\eqs
has finitely many zeros in the strip $\mathcal{S}_\eta$, by analyticity and boundedness of $\Delta_{\A^{\rho_0}}(\nu)$. We define the \textbf{crossing number} of $\A$,  $\mathrm{cross}(\A)$, to be the net number of roots (counted with multiplicity) which cross the imaginary axis from left to right as $\rho$ increases from $-\infty$ to $+\infty$. More precisely, fix any $\rho_j \in \text{NH}(\A)$ and let $\left(\nu_{j,l}\right)_{l=1}^{k_j}$ denote the roots of  $d_{\rho_j}(\nu)$ on the imaginary axis, $\Re(\nu_{j,l})=0$. We list multiple roots repeatedly according to their multiplicity.  Let $M_j$ denote the sum of their multiplicities. For $\rho$ near $\rho_j$, with $\pm(\rho-\rho_j)>0$, this equation has exactly $M_j$ roots (counting multiplicity) near the imaginary axis, $M_j^{L_\pm}$ with $\Re \nu <0$ and $M_j^{R_\pm}$ with $\Re \nu >0$, and $M_j=M_j^{L_\pm}+M_j^{R_\pm}$. We define
\bqs
\mathrm{cross}(\A)=\sum_{j=1}^m\left(M_j^{R_+}-M_j^{R_-} \right).
\eqs

For $\A\in\cP^1$, we say that a crossing $\rho_0$ is \textbf{simple} if there is precisely one simple root of $d_{\rho_j}(\nu_*)$ located on the imaginary axis, and if this root crosses the imaginary axis with non-vanishing speed as $\rho$ passes through $\rho_j$. Note that for these simple crossings, we can locally continue the root  $\nu_*\in i\R$ as a $C^1$-function of $\rho$ as $\nu(\rho)$. We refer to this root as the \textbf{crossing root}. Non-vanishing speed of crossing can then be  expressed as $ \Re\left( \dot \nu(\rho_0)\right) \neq 0$. 

Next, suppose that $\A\in\cP^1$ has only simple crossings $\rho_j\in\text{NH}(\A)$. In this case we let $\nu_j(\rho)$ be the complex-valued crossing-value defined near $\rho_j$ such that $\nu_j(\rho)$ is a root of $d_{\rho}$ and $\Re(\nu_j(\rho_j))=0$. In this case, the crossing number is explicitly given through
\bqq
\label{crossNum}
\mathrm{cross}(\A)=\sum_{j=1}^m\text{sign}\left( \Re\left( \dot \nu_j(\rho_j)\right)\right).
\eqq

The following result shows that the set of paths with only simple crossings is dense in $\cP$.

\begin{lem}\label{lemdensity}
Let $\A\in\cP$ be such that $\text{NH}(\A)$ is a finite set. Then given $\epsilon>0$, there exists $\widetilde{\A}\in\cP^1$ such that:
\begin{itemize}
\item[(i)] $\widetilde{\A}^\pm=\A^\pm$;
\item[(ii)] $|\widetilde{\A}^\rho-\A^\rho|<\epsilon$ for all $\rho\in\R$; and
\item[(iii)] $\widetilde{\A}$ has only simple crossings.
\end{itemize}
\end{lem}

This lemma is proved in the following section.

\begin{rmk}
If $\epsilon$ is small enough in Lemma \ref{lemdensity}, then one has
\bqs
\mathrm{cross}(\A)=\mathrm{cross}(\widetilde{\A}).
\eqs
\end{rmk}

\subsection{Proof of Lemma \ref{lemdensity}}

The proof follows \cite{mallet-paret:99} with some appropriate modifications. 

We start by introducing submanifolds of $\M_n(\C)$. For $0\leq k \leq n$ we define the sets $\mG_k\subset \M_n(\C)$ and $\mH\subset \M_n(\C)\times\M_n(\C)$ by
\begin{align*}
\mG_k&=\left\{M\in\M_n(\C)~|~ \text{rank}(M)=k \right\},\\
\mH&=\left\{(M_1,M_2)\in\M_n(\C)\times\M_n(\C) ~|~ \text{rank}(M_1)=n-1, \right.\\
&~\quad \left.M_2\text{ is invertible, and } \text{rank}(M_1M_2^{-1}M_1)=n-2 \right\}.
\end{align*}

The sets $\mG_k$ and $\mH$ are analytic submanifolds of $ \M_n(\C)$ and $\M_n(\C)\times\M_n(\C)$ respectively, of complex dimension
\bqq
\label{dim}
\dim_\C \mG_k=n^2-(n-k)^2, \quad \dim_\C \mH=2n^2-2;
\eqq
see \cite{mallet-paret:99}.
We also consider the following maps
\begin{align*}
\F,\G&:\left(L^1_\eta(\R,\M_n(\C))\times \ell^1_\eta\left( \M_n(\C)\right)\right)\times\R\rightarrow \M_n(\C)\\
\F\times\G&: \left(L^1_\eta(\R,\M_n(\C))\times\ell^1_\eta\left(\M_n(\C) \right)\right)\times\R\rightarrow \M_n(\C)\times \M_n(\C)\\ 
\mathcal{D}&: \left(L^1_\eta(\R,\M_n(\C))\times\ell^1_\eta \left(\M_n(\C) \right)\right) \times T \rightarrow \M_n(\C)\times \M_n(\C)
\end{align*}
given by
\begin{subequations}
\label{mapsFGD}
\begin{align}
\F(\A,\ell)&=i\ell \mathbb{I}_n-\widehat{\K}(i\ell)-\sum_{j\in\J}A_je^{-i\ell \xi_j},\\
\G(\A,\ell)&=\mathbb{I}_n-\widehat{\K}'(i\ell)+\sum_{j\in\J}\xi_j A_je^{-i\ell \xi_j},\\
(\F\times\G)(\A,\ell)&=\left(\F(\A,\ell),\G(\A,\ell)\right),\\
\mathcal{D}(\A,\ell_1,\ell_2)&=\left(\F(\A,\ell_1),\F(\A,\ell_2)\right),
\end{align}
\end{subequations}

where $\A=\left(\K,\left(A_j\right)_{j\in\J}\right)\in L^1_\eta(\R,\M_n(\C))\times \ell^1_\eta \left(\M_n(\C) \right)$ and $T$ is the set
\bqs
T=\left\{(\ell_1,\ell_2)\in\R^2~|~ \ell_1<\ell_2 \right\}.
\eqs

\begin{prop}

Suppose that $\A=\left(\K,\left(A_j\right)_{j\in\J}\right)\in L^1_\eta(\R,\M_n(\C))\times \ell^1_\eta \left(\M_n(\C) \right)$ satisfies the conditions
\bqq
\label{conditionsFGD}
\begin{matrix}
(i) & \F(\A,\ell)\notin \mG_k, & 0\leq k \leq n-2, & \ell \in \R \\
(ii) & (\F\times\G)(\A,\ell)\notin \mG_{n-1}\times\mG_k, &  0\leq k \leq n-1, & \ell \in \R \\
(iii) & (\F\times\G)(\A,\ell)\notin\mH, & \ell\in\R & \\
(iv) & \mathcal{D}(\A,\ell_1,\ell_2)\notin \mG_{n-1}\times\mG_{n-1},& (\ell_1,\ell_2)\in T &
\end{matrix}
\eqq
for all ranges of $k,\ell,\ell_1$ and $\ell_2$. Then the constant coefficient system \eqref{eqNLconstant} has at most one $\ell\in\R$ such that $\nu=i\ell$ is a root of the characteristic equation $\det \Delta _{\A}(\nu)=0$, and the root $\nu$ is simple. 
\end{prop}

\begin{Proof}
We first note that $\F(\A,\ell)=\Delta_{\A}(i\ell)$ as defined in \eqref{matDconstant} and that $\G(\A,\ell)=-i\Delta_{\A}'(i\ell)$. Therefore, condition $(i)$ implies that $\text{rank}(\Delta_\A(\nu))=n-1$ for all roots $\nu=i\ell$. Condition $(ii)$ ensures that $\Delta_\A'(\nu)$ is invertible for such $\nu$. Condition $(iii)$ implies that the rank of $\Delta_\A(\nu)\Delta_\A'(\nu)^{-1}\Delta_\A(\nu)$ is $n-1$ for such $\nu$. Hypothesis \ref{bound} ensures  the existence of $\eta_0>0$ such that $\eta-\eta_0>0$ and  $f(\nu)=\Delta_\A(\nu)$ is a holomorphic function in a neighborhood of $i\ell \in \mathcal{S}_{\eta-\eta_0}=\left\{\nu\in\C~|~|\Re(\nu)|<\eta-\eta_0 \right\}$ that satisfies:
\begin{itemize}
\item $\text{rank}(f(i\ell))=n-1$
\item $f'(i\ell)$ is invertible
\item $\text{rank}(f(i\ell)f'(i\ell)^{-1}f(i\ell))=n-1$.
\end{itemize}
As a consequence, $g(\nu)=\det f(\nu)$ has a simple root at $\nu=i\ell$ \cite{mallet-paret:99} and $\nu=i\ell$ is a simple root of the characteristic equation $\det \Delta_\A(\nu)=0$. Finally, the last condition $(iv)$ ensures that there is at most one value $\ell\in\R$ for which $\det \Delta_\A(i\ell)=0$ which concludes the proof.
\end{Proof}

\begin{prop}\label{propsurjective}
The maps $\F$ and $\F\times\G$ have surjective derivative with respect to the first argument $\A$ at each point $(\A,\ell)\in L^1(\R,\M_n(\C))\times \ell^1 \left(\M_n(\C) \right)\times\R$. Moreover, if $\xi_j/\xi_k$ is irrational for some $j<k$, then the derivative of the map $\mathcal{D}$  with respect to the first argument $\A$ is surjective at each $(\A,\ell)\in L^1(\R,\M_n(\C))\times \ell^1 \left(\M_n(\C) \right)\times T$.
\end{prop}

\begin{Proof}
From their respective definition, one sees immediately that the derivative of $\F$ with respect to $A_1\in\M_n(\C)$ is $-\mathbb{I}_n$ and that the derivative with respect to $(A_1,A_2)\in\M_n(\C)\times\M_n(\C)$ is given by the matrix
\bqs
-\left(\begin{matrix} \mathbb{I}_n & e^{-i\ell \xi_2} \mathbb{I}_n\\
0_n & -\xi_2  e^{-i\ell \xi_2} \mathbb{I}_n
\end{matrix} \right)
\eqs
which is an isomorphism on $\M_n(\C)\times \M_n(\C)$; in particular, the derivative of both maps is onto.

We fix $(\ell_1,\ell_2)\in T$. Then at least one of the quantities $(\ell_1-\ell_2)\xi_j$ or $(\ell_1-\ell_2)\xi_k$ is irrational. Suppose now that $(\ell_1-\ell_2)\xi_j$ is irrational. Then the derivative of $\mathcal{D}$ with respect to $(A_1,A_j)$ is given by
\bqs
-\left(\begin{matrix} \mathbb{I}_n & e^{-i\ell_1 \xi_j} \mathbb{I}_n\\
\mathbb{I}_n & e^{-i\ell_2 \xi_j} \mathbb{I}_n
\end{matrix} \right)
\eqs
which is an isomorphism.
\end{Proof}

\begin{rmk}
Note that we can always assume that $\xi_j/\xi_k$ is irrational for some $j<k$. If it is not the case, we can enlarge $\J$ to $\J\cup \left\{\xi_* \right\}$ with an additional constant coefficient $A_{*}=0$ in  \eqref{eqNL} so that $\xi_*/\xi_k$ is irrational for some $k\in\J$.
\end{rmk}

In order to complete the proof of Lemma \ref{lemdensity}, we will use the notion of transversality for smooth maps defined in manifolds. We say that a smooth map $f:\X\rightarrow\Y$ from two manifolds is transverse to a submanifold $\mathcal{Z}\subset \Y$ on a subset $\mathcal{S}\subset\X$ if
\bqs
\text{rg}(Df(x))+T_{f(x)}\mathcal{Z}=T_{f(x)}\Y \quad \text{ whenever } x\in\mathcal{S} \text{ and } f(x)\in \mathcal{Z}
\eqs
where $T_pM$ denotes the tangent space of $M$ at a point $p$. 

\begin{thm}[Transversality Density Theorem]\label{thmtransversality} Let $\mathcal{V},\X,\Y$ be $\cC^r$ manifolds, $\Psi:\mathcal{V}\rightarrow \cC^r(\X,\Y)$ a representation and $\mathcal{Z}\subset\Y$ a submanifold and $\text{ev}_\Psi:\mathcal{V}\times \X\rightarrow \Y$ the evaluation map. Assume that:
\begin{enumerate}
\item $\X$ has finite dimension $N$ and $\mathcal{Z}$ has finite codimension $Q$ in $\Y$;
\item $\mathcal{V}$ and $\X$ are second countable;
\item $r>\max(0,N-Q)$;
\item $\text{ev}_\Psi$ is transverse to $\mathcal{Z}$.
\end{enumerate}
Then the set $\left\{ V\in \mathcal{V}~|~ \Psi_V \text{ is transverse to } \mathcal{Z} \right\}$ is residual in $\mathcal{V}$.
\end{thm}

The proof of this theorem can be found in \cite{abraham-robbin:70}.

\begin{prop}\label{propresidual}
There exists a residual (and hence dense) subset  of $\cP^1$ such that for any $\A$ in this subset, all conditions \eqref{conditionsFGD} are satisfied.
\end{prop}

\begin{Proof}
The idea is to apply the Transversality Density Theorem \ref{thmtransversality} to exhibit a residual subset of $\cP^1$ such that all the maps $\F(\A^\rho,\ell),(\F\times\G)(\A^\rho,\ell)$ and $\mathcal{D}(\A^\rho,\ell_1,\ell_2)$ are transverse to the manifolds appearing in \eqref{conditionsFGD} on $(\rho,\ell)\in\R^2$ and $(\rho,\ell_1,\ell_2)\in\R^2$ respectively. For simplicity we only detail the proof for $\F$, the two other cases being similar.

We apply Theorem \ref{thmtransversality} with manifolds $\mathcal{V}=\cP^1$, $\X=\R^2$ and $\Y=\M_n(\C)$ and submanifold $\mathcal{Z}=\mG_k$ with $0\leq k \leq n-2$. So for any $\A\in\cP^1$ we define $\Psi_\A:\R^2\rightarrow \M_n(\C)$ by
\bqs
\Psi_\A(\rho,\ell)=\F(\A^\rho,\ell),
\eqs
and the evaluation map is simply given by $\text{ev}_\Psi:\cP^1\times\R^2\rightarrow \M_n(\C)$
\bqs
\text{ev}_\Psi(\A,\rho,\ell)=\F(\A^\rho,\ell).
\eqs
We thus have $r=1$, $N=2$ and $Q=2(n-k)^2$ (the real codimension of $\mG_k$). This implies that the third condition of Theorem \ref{thmtransversality} is satisfied for all $0\leq k \leq n-2$. Proposition \ref{propsurjective} ensures that the required transversality hypothesis of the evaluation map is fulfilled.

We can then conclude that there exists a residual subset (and hence dense) of $\cP^1$ such that for any $\A$ in this subset the composed map $\F(\A^\rho,\ell)$ is transverse to the manifolds appearing in \eqref{conditionsFGD}.
\end{Proof}

\begin{Proof}[of Lemma  \ref{lemdensity}]
We are now ready to prove Lemma \ref{lemdensity}. Let $\A\in \cP$ such that $NH(\A)$ is a finite set. By Proposition \ref{propresidual}, we may assume that the family $\A$ in the statement of Lemma \ref{lemdensity} is such that all four conditions \eqref{conditionsFGD} hold for $\A^\rho$ for each $\rho\in\R$. Thus for each such $\A^\rho$, the constant coefficient equation \eqref{eqNLconstant} has at most one $\ell\in\R$ such that $\nu=i\ell$ is an root and $i\ell$ is a simple root of the characteristic equation $\det \Delta _{\A^\rho}(\nu)=0$. It is then enough to perturb $\A$ to a nearby $\widetilde{\A}\in \cP^1$ with the same endpoints $\widetilde{A}^\pm=\A^\pm$ such that, by Sard's Theorem, all the roots of the corresponding family of equations \eqref{eqNLconstant} cross the imaginary axis transversely with $\rho$, that is, $\widetilde{A}$ has only simple crossings.
\end{Proof}

\subsection{Proof of Theorem \ref{thmSFI}}
We first introduce the map $\Sigma_\gamma:L^1_\eta(\R,\M_n(\C))\times \ell^1_\eta \left(\M_n(\C) \right)\rightarrow L^1_\eta(\R,\M_n(\C))\times \ell^1_\eta \left(\M_n(\C) \right)$, defined for each $\gamma\in\R$ by
\bqs
\Sigma_\gamma \cdot \A^0=\Sigma_\gamma \cdot \left( \K^0, \left(A_j^0 \right)_{j\in\J}\right):= \left( \K^0_\gamma, \left(A_{j,\gamma}^0 \right)_{j\in\J}\right),
\eqs
where
\bqs
\K_\gamma^0(\zeta)=\K^0(\zeta)e^{\gamma \zeta},\quad \forall\, \zeta \in \R, \quad
A_{1,\gamma}^0= A_1^0+\gamma, \quad
A_{j^,\gamma}^0=A_j^0e^{\gamma \xi_j},\quad \forall\, j\neq 1.
\eqs
This transformation $\Sigma_\gamma$ arises from a change of variables $V(\xi)=e^{\gamma\xi}U(\xi)$ in \eqref{eqNLconstant} with constant coefficient $\A^0= \left( \K^0, \left(A_j^0 \right)_{j\in\J}\right)$. One can then easily check that
\bqs
\Delta_{\Sigma_\gamma \cdot \A^0}(\nu)=\Delta_{\A^0}(\nu-\gamma), \quad \nu\in\C,
\eqs
so that $\Sigma_\gamma$ shifts all eigenvalues to the right by an amount of $\gamma$.

\begin{prop}\label{propSC}
Suppose that $\nu=i\ell$, with $\ell\in\R$, is a simple root of the characteristic equation \eqref{eqcharacteristic} associated to $\A^0$, and suppose that there are no other roots with $\Re \lambda =0$. Then for $\gamma\in\R$, $0<|\gamma| <\eta $ sufficiently small, we have that
\bqq
\iota(\Sigma_{-\gamma}\cdot\A^0,\Sigma_{\gamma}\cdot\A^0)=-\mathrm{sign}(\gamma).
\eqq
\end{prop}

\begin{Proof}
With $\A^0=\left( \K^0, \left(A_j^0 \right)_{j\in\J}\right)$, we make the change of variable $V(\xi)=W_\gamma(\xi)U(\xi)$, with $W_\gamma(\xi)=e^{\gamma\sqrt{\xi^2+1}}$, in equation \eqref{eqNLconstant}. This leads to a nonautonomous equation of the form 
\bqs
V'(\xi)=\cN[\A_\gamma(\xi)]\cdot V(\xi)
\eqs
which is asymptotically hyperbolic with limiting operators $\A_\gamma^\pm=\Sigma_{\pm\gamma}\cdot \A^0$. It is easy to check that the exponential localization of $\K^0$ in $L^1_\eta(\R,\M_n(\C))$ and $\left(A_{j,\gamma}^0 \right)_{j\in\J}$ in $\ell^1_\eta(\R,\M_n(\C))$ ensures that $\T_{\A_\gamma}$ satisfies Hypotheses \ref{hypK} and \ref{hypA}, and thus is Fredholm for $0<\gamma<\eta$. Similarly, we make the change of variable $V(\xi)=W_{-\gamma}(\xi)U(\xi)$ in the adjoint equation
\bqs
U'(\xi)=\cN[\A^0]^*U(\xi)
\eqs
which results in the nonautonomous equation
\bqs
V'(\xi)=\cN[\A_\gamma(\xi)]^*\cdot V(\xi).
\eqs
Without loss of generality, we suppose that $\gamma>0$ and is small enough so that $\nu=i\ell$ is the only root of $\det(\Delta_{\A^0})=0$ in the strip $|\Re(\nu)|\leq\gamma<\eta$. Suppose that $V$ is a nonzero element of the kernel of $\T_{\A_\gamma}$, then $V$ is bounded and $U=W_{-\gamma}V$ is also a bounded solution of $U'(\xi)=\cN[\A^0]\cdot U(\xi)$, hence $U(\xi)=e^{i\ell \xi}p$ for some nonzero vector $p\in\C^n$. Indeed, as $\nu=i\ell$ is a simple root of $\det(\Delta_{\A^0})=0$, there exists $p\in\C^n$, such that $p$ belongs to the kernel of $\Delta_{\A^0}$ and thus $e^{i\ell\xi}p$ is in the kernel of $\T_{\A^0}$ using Fourier transform. But, $V(\xi)=W_{\gamma}(\xi)e^{i\ell \xi}p$ is now unbounded which is a contradiction, and so $\text{ker}~\T_{\A_\gamma}=\left\{0 \right\}$. Applying a similar argument to the adjoint equation, one sees that $V(\xi)=W_{-\gamma}(\xi)e{-i\ell\xi}q$, $q\in\C^n$, is the one-dimensional span of $\text{ker}~\T_{\A_\gamma}^*$. Thus, applying Theorem \ref{thmFredholm}, we have
\bqs
\iota(\Sigma_{-\gamma}\cdot\A^0,\Sigma_{\gamma}\cdot\A^0)=\text{ind}~\T_{\A_\gamma}=\dim\text{ker}~\T_{\A_\gamma} - \dim \text{ker}~\T_{\A_\gamma}^*=-1.
\eqs
\end{Proof}

The following proposition shows that without loss of generality we may assume that roots of the characteristic equation cross the imaginary axis by means of a rigid shift of the spectrum with the operator $\Sigma_\gamma$. 

\begin{prop}\label{propShift}
Let $\A\in\cP^1$ be such that $NH(\A)$ is a finite set and has only simple crossings. Then there exists $\widetilde{\A}\in\cP^1$ such that:
\begin{itemize}
\item[(i)] $\A^\pm=\widetilde\A^\pm$;
\item[(ii)] $NH(\A)=NH(\widetilde\A)$;
\item[(iii)] for each $\rho_j\in NH(\A)$, we have $\Re(\dot\nu_j(\rho_j))=\Re(\dot{\widetilde\nu}_j(\rho_j))$, with $\widetilde{\nu}_j$ corresponding to $\widetilde{\A}$;
\item[(iv)] $\widetilde\A$ has only simple crossings.
\end{itemize}
In addition, the family $\widetilde{A}$ has the form
\bqq
\label{familywA}
\widetilde{\A}^\rho=\Sigma_{\gamma_j(\rho-\rho_j)}\cdot \A^{\rho_j}, \quad \gamma_j:=\Re(\dot\nu_j(\rho_j)),
\eqq
for $\rho$ in a neighborhood of each $\rho_j$.
\end{prop}

We omit the proof of this result, as it is identical to that in \cite{mallet-paret:99}.

\begin{Proof}[of Theorem \ref{thmSFI}] 
Let $\A\in\cC\left(\R,L^1_\eta(\R,\M_n(\C))\times \ell^1_\eta \left(\M_n(\C) \right)\right)$ be a one-parameter family as in the statement of Theorem \ref{thmSFI}. Without loss, by Lemma \ref{lemdensity}, we may assume that $\A$ has only simple crossings. Let $\widetilde\A\in\cC^1\left(\R,L^1_\eta(\R,\M_n(\C))\times \ell^1_\eta \left(\M_n(\C) \right)\right)$ as in statement of Proposition \ref{propShift}. Then for any sufficiently small $\epsilon>0$, using the Corollary \ref{cor:cocycle}, we have that
\bqs
\iota(\A^-,\A^+)=\iota(\A^-,\widetilde{\A}^{\rho_1-\epsilon})+\sum_{j=1}^{m-1}\iota(\widetilde\A^{\rho_j+\epsilon},\widetilde\A^{\rho_{j+1}-\epsilon})+\sum_{j=1}^m\iota(\widetilde\A^{\rho_j-\epsilon},\widetilde\A^{\rho_{j}+\epsilon})+\iota(\widetilde\A^{\rho_m+\epsilon},\A^+).
\eqs
For each $\rho$ in the intervals: $[\rho_j+\epsilon,\rho_{j+1}-\epsilon]$, $1\leq j \leq m-1$, $(-\infty,\rho_1-\epsilon]$ and $[\rho_n+\epsilon,+\infty)$, equation \eqref{eqNLconstant} is hyperbolic, and one concludes that
\bqs
\iota(\A^-,\widetilde{\A}^{\rho_1-\epsilon})=\sum_{j=1}^{m-1}\iota(\widetilde\A^{\rho_j+\epsilon},\widetilde\A^{\rho_{j+1}-\epsilon})=\iota(\widetilde\A^{\rho_m+\epsilon},\A^+)=0.
\eqs 
On each interval $[\rho_j-\epsilon,\rho_{j}+\epsilon]$, $1\leq j \leq m$, we have a simple crossing and we can apply the result of Proposition \ref{propSC}:
\bqs
\sum_{j=1}^m\iota(\widetilde\A^{\rho_j-\epsilon},\widetilde\A^{\rho_{j}+\epsilon})=-\sum_{j=1}^m\text{sign}\left( \Re\left( \dot \nu_j(\rho_j)\right)\right).
\eqs
This implies that $\iota(\A^-,\widetilde{\A}^{\rho_1-\epsilon})=-\text{cross}(\A)$ which concludes the proof.
\end{Proof}

\subsection{Exponentially weighted spaces}

We now give a first application of Theorem \ref{thmFlow} to operators posed on exponentially weighted spaces. Assume that $\A^0=\left( \K^0, \left(A_j ^0\right)_{j\in\J}\right)\in\cP$ is a constant coefficient operator and consider the associated operator $\T_{\A^0}=\frac{d}{d\xi}-\cN[\A^0]$ on the space $\widetilde{L}^2_\eta(\R,\C^n)$ with norm
\bqs
\| U \|_{\widetilde{L}^2_\eta}=\|U(~\cdot~)e^{\gamma~ \cdot ~ }\|_{L^2(\R,\C^n)}.
\eqs 
Using the isomorphism 
\bqs
\widetilde{L}^2_\gamma(\R,\C^n) \longrightarrow L^2(\R,\C^n), \quad U(\xi)\longmapsto U(\xi)e^{\gamma \xi},
\eqs
the operator $\T_{\A^0}$ for $U$ on $\widetilde{L}^2_\gamma(\R,\C^n)$ is readily seen to be conjugate to $\T_{\A^0}^\gamma=\frac{d}{d\xi}-\cN[\Sigma_\gamma\cdot\A^0]$ for $V$ on $L^2(\R,\C^n)$. We conclude that $\T_{\A^0}^\gamma$ is Fredholm for $\gamma$ in open subsets of the real line. When $\A^0$ has only finitely many simple crossings, we can consider the family of operators  $\T_{\A^0}^\gamma$ with $\gamma$ close to zero. More generally, we introduce a two-sided family of weights via
\bqs
\| U \|_{\gamma_-,\gamma_+}=\| U \chi_+ \|_{\widetilde{L}^2_{\gamma_+}} + \| U \chi_- \|_{\widetilde{L}^2_{\gamma_-}}
\eqs
where
\bqs
\chi_\pm(\xi)=\left\{
\begin{matrix}
1 & \pm \xi >0 \\
0 & \text{otherwise}.
\end{matrix}
 \right.
\eqs
The operator $\T_{\A^0}$ on $L^2_{\gamma_-,\gamma_+}$ is conjugate to an operator $\T_{\A^0}^{\gamma_-,\gamma_+}$ on $L^2$ whose coefficients are $\Sigma_{\gamma_+}\cdot\A^{0}$ for $\xi>0$ and $\Sigma_{\gamma_-}\cdot \A^{0}$ for $\xi<0$. The following corollary is a direct consequence of the above discussion and Theorem \ref{thmFlow}.

\begin{cor}\label{corweight}
Suppose that $\nu=i\ell$, with $\ell\in\R$, is a root of the characteristic equation associated to $\A^0$ of multiplicity $N$, and suppose that there are no other roots with $\Re \lambda =0$. Then, the operator $\T_{\A^0}^{\gamma_-,\gamma_+}$ is Fredholm for all $\gamma_\pm$ close to zero with $\gamma_-\gamma_+\neq 0$ and for $\gamma\in\R$, $\gamma\neq 0$ sufficiently small, we have that
\bqq
\iota(\Sigma_{-\gamma}\cdot\A^0,\Sigma_{\gamma}\cdot\A^0)=\mathrm{ind}~\T_{\A^0}^{-\gamma,\gamma}=-\mathrm{sign}(\gamma)N.
\eqq
\end{cor}

\section{Applications}\label{application}

We give two applications of our main result. We first consider the effect of small inhomogeneities in nonlocal conservation laws. We then show how our results can be used to study edge bifurcations for nonlocal eigenvalue problems, replacing Gap Lemma constructions with Lyapunov-Schmidt and far-field matching constructions.

\subsection{Localized source terms in nonlocal conservation laws}\label{secNCL}

Consider the nonlocal conservation laws 
\bqq\label{e:nlcl}
U_t =\left( \K \ast F(U) +  G(U) \right)_x,\quad U\in \R^n,\ x\in\R,
\eqq
with appropriate conditions on convolution kernel $\K$, and fluxes $F,G$. Nonlocal conservation laws arise in a variety of applications and pose a number of analytic challenges; see \cite{du-etal:12} for a recent discussion and references. 

In the absence of the nonlocal, dispersive term $\K\ast F$, the system of conservation laws is well known to develop discontinuities in finite time which are referred to as shocks. Shocks can usually be classified according to ingoing and outgoing characteristics. In the presence of viscosity, shocks are smooth traveling waves, and characteristic speeds can be characterized via the group velocities of neutral modes in the linearization. In our case, the linearization at a constant state 
\[
V_t =\left( \K \ast dF_U(0) +  dG_U(0) \right)V_x,\quad V\in \R^n,\ x\in\R,
\]
can be readily solved via Fourier transform, with dispersion relation \[
d(\lambda,i\ell)=\mathrm{det}\,\left( i\ell\widehat{\K}(i\ell) dF_U(0) + i\ell dG_U(0)-\lambda\mathbb{I}_n\right).
\]
We find an eigenvalue $\lambda=0$ with multiplicity $n$. Assuming that $\widehat{\K}(i\ell) dF_U(0) +  dG_U(0)$ possesses real, distinct eigenvalues $-c_j$, we obtain expansions $\lambda_j(i\ell)=-c_j\ell+\rmO(\ell^2)$, so that the negative eigenvalues $c_j$ naturally denote speeds of transport in different components of the system. As with viscous approximations to local conservation laws, instabilities can enter for finite wavenumber $\ell$ for non-scalar diffusion, so that we will need an extra condition on the nonlocal part that guarantees stability of the homogeneous solution. 

Rather than studying existence of large-amplitude shock profiles, we focus here on a perturbation result, exploiting the linear Fredholm theory developed in the previous sections. It will be clear from the techniques employed here and in the subsequent section that our results can be used to develop a spectral theory for large amplitude shock profiles in the spirit of \cite{hz}. Our results parallel the results in \cite{sandstedescheel:08}, where viscous regularization of conservation laws were analyzed. Roughly speaking, our results show that at small amplitude, nonlocal, dispersive terms act in a completely analogous fashion to viscous regularizing terms. 

Our analysis considers spatially localized source terms of the nonlocal conservation law \eqref{e:nlcl},
\bqq
\label{eq:nvcl}
U_t =\left( \K \ast F(U) +  G(U) \right)_x + \epsilon H(x,U,U_x),\quad U\in \R^n
\eqq
for a kernel $\K \in L^1_{\eta_0} (\R,\M_n(\R))$, with fixed $\eta_0>0$, and a smooth hyperbolic flux $g$ with
\begin{subequations}
\label{hypflux}
\begin{align}
\det (dG_U(0)) & \neq 0 \label{hypflux1} \\
\sigma \left(dG_U(0) +\widehat{\K}(0)dF_U(0) \right)&= \left\{-c_1>-c_2>\dots>-c_n \right\} \label{hypflux2} \\
\det \left( dG_U(0)+\widehat{\K}(i\ell) dF_U(0) \right) &\neq 0, \quad \forall\, \ell \in \R, \ell\neq 0 \label{hypflux3}\\
\widehat{\K}(\nu) dF_U(0),dG_U(0)&\in \mathcal{S}_n(\R)=\left\{ M \in \M_n(\R) ~|~ M=M^t \right\} \quad \forall\, \nu \in \C
\end{align}
\end{subequations}
and a smooth, spatially localized, source term $H$ so that there exist constant $C,\delta>0$ such that
\bqq
\label{eq:source}
\| H(x,U,V) \| \leq C e^{-\delta |x|}
\eqq
for all $x\in\R$ and all $(U,V)$ near zero in $\R^n\times\R^n$. 

Here, the first condition guarantees that steady-states are solutions to ODEs, hence smooth; the second condition enforces strict hyperbolicity of the nonlocal linear part, the third condition guarantees that zero is not in the essential spectrum of the linearization for any nonzero wavenumber. The last condition refers to the usual requirement of symmetric fluxes. 

We look for small bounded solutions of the nonlocal equation
\bqq
\label{eq:nvclstat}
0 =\left( \K \ast F(U) +  G(U) \right)_x + \epsilon H(x,U, U_x).
\eqq

Contrary to hyperbolic conservation laws where the viscous term is typically $B U_{xx}$ with a positive definite, symmetric viscosity matrix $B$, we cannot use spatial dynamics techniques for \eqref{eq:nvclstat} because of the nonlocal term $\K \ast F(U)$. Instead, following \cite{sandstedescheel:08}, we will use an approach based only on functional analysis and Lyapunov-Schmidt reduction, thus exploiting the Fredholm and spectral flow properties developed in the previous sections. The key point of our approach is the linearization of equation \eqref{eq:nvclstat} at the solution $U=0$ and $\epsilon = 0$
\bqq
\label{eq:linear}
\cl U = \K_x \ast \left( dF_U(0) U \right)+dG_U(0) U_x.
\eqq
The adjoint $\cl^*$ of \eqref{eq:linear} is given by
\bqq
\label{eq:linear2}
\cl^* U = - dF_U(0)^t \K_-^t \ast U_x -dG_U(0)^t U_x
\eqq
where $\K_-^t(x)=\K^t(-x)$. Assuming that $dG_U(0)$ is invertible, we can associate the operator 
\bqq
\label{eq:linearmod}
\widetilde{\cl} U = U_x+ dG_U(0)^{-1}\K_x \ast \left( dF_U(0) U \right)
\eqq
which is of the form of a constant operator studied in Section \ref{fredholm} as $\K_x\in L^1 (\R,\M_n(\R))$. Both $\cl$ and $\widetilde{\cl}$ can be viewed as unbounded linear operators on $L^2(\R,\R^n)$ but also can be considered as unbounded operators on $L^2_\eta(\R,\R^n)$ for $0<\eta<\eta_0$ as $\K \in L^1_{\eta_0} (\R,\M_n(\R))$ with norm
\bqs
\| U \|_{L^2_\eta(\R,\R^n)}=\| U(x) e^{\eta |x| }\|_{L^2(\R,\R^n)}.
\eqs

\begin{lem}\label{lem1}
Assume that $c_j \neq 0$ for all $j$, then there is an $\eta_*>0$ with the following property. For each fixed $\eta$ with $0<\eta<\eta_*$, the operator $\cl$ defined on $L^2_\eta(\R,\R^n)$ is Fredholm with index $-n$ and has trivial null space.
\end{lem}

\begin{Proof}
The characteristic equation associated to the linearized system \eqref{eq:linearmod} is
\bqq
\label{eq:charalin}
0=\det(\nu \mathbb{I}_n+ \nu dG_U(0)^{-1}\widehat{\K}(\nu)dF_U(0))=\nu^n \det(dG_U(0)^{-1}) \det \left( dG_U(0)+\widehat{\K}(\nu) dF_U(0) \right),
\eqq
so that $\nu=0$ is an root with multiplicity $n$, and all other roots have nonzero real part due to \eqref{hypflux3}. We can apply Corollary \ref{corweight} and find that the Fredholm index of $\widetilde{\cl}$ and thus of $\cl$ on $L^2_\eta(\R,\R^n)$ is equal to $-n$ as claimed. Since \eqref{eq:linear} is translation invariant, we can use Fourier transform to analyze the kernel.  Any function $U$ in the kernel of $\cl$ satisfies
\bqs
0=i\ell \left( \widehat{\K}(i\ell)dF_U(0)+ dG_U(0)\right) \widehat{U}(\ell).
\eqs
As $U \in L^2_\eta(\R,\R^n)$, $\widehat{U}(\ell)$ is a bounded analytic function in the strip $\mathcal{S}_\eta$, and thus $\widehat{U}(\ell)=0$ for all $\ell \in \R$. This proves that the kernel of $\cl$ in the exponentially weighted space is trivial.
\end{Proof}

Lemma \ref{lem1} implies that the kernel of the $L^2$-adjoint $\cl^*$ of $\cl$ considered on $L^2_{-\eta}(\R,\R^n)$ is $n$-dimensional and thus spanned by the constants $e_j$ for $j=1,\dots,n$ where $e_j$ form an orthonormal basis of $\R^n$ such that
\bqs
\left( dG_U(0)+\widehat{\K}(0)dF_U(0)\right)e_j=-c_je_j.
\eqs

To find shock-like transition layers, caused by the inhomogeneity $h$ for small $\epsilon$, we make the following ansatz
\bqq
\label{ansatzappli}
U(x)=\sum_{j=1}^na_je_j\chi_+(x)+\sum_{j=1}^n b_j e_j \chi_-(x)+W(x),
\eqq
where $a_j,b_j\in\R$ and $W\in L^2_\eta(\R,\R^n)$, and $\chi_\pm(x)=(1+\tanh(\pm x))/2$. Substituting the ansatz into \eqref{eq:nvcl}, we obtain an equation of the form
\bqq
\label{eq:functionF}
\F(a,b,W;\epsilon)=0, \quad \F(~\cdot~;\epsilon):\R^n\times\R^n\times \mathcal{D}(\cl) \longrightarrow L^2_\eta(\R,\R^n)
\eqq
for $a=(a_j),b=(b_j)$. For small enough $\eta$, the map $\F$ is smooth and the its linearization at $(a,b,W)=0$ is given by
\bqs
\F_W(0;0)=\cl, \quad \F_{a_j}(0,0)=\K_x\ast (dF_U(0)e_j \chi_+ )+dG_U(0)e_j\chi', \quad \F_{b_j}(0,0)=\K_x\ast (dF_U(0)e_j \chi_- )-dG_U(0)e_j\chi'
\eqs
where $\F_a(0,0)$ and $\F_b(0,0)$ lie in $L^2_\eta(\R,\R^n)$.

\begin{lem}
Under the hypotheses of Lemma \ref{lem1}, the operator
\bqs
\F_{a,W}(0;0):\R^n\times L^2_\eta(\R,\R^n) \longrightarrow L^2_\eta(\R,\R^n), \quad (a,W) \longmapsto \F_a(0;0)a+\F_W(0;0)W
\eqs
is invertible.
\end{lem}

\begin{Proof}
We first note that the $n$ partial derivatives with respect to $a_j$ are linearly independent. To see this, we integrate $\F_{a_j}(0,0)$ over the real line to find
\begin{align*}
\int_\R \F_{a_j}(0,0)dx &=\int_\R \K_x\ast \left(dF_U(0)e_j \chi_+\right) dx +dG_U(0)e_j \\
&=\int_\R \K \ast \left(dF_U(0)e_j \chi'_+\right)dx +dG_U(0)e_j \\
&= \widehat{\K}(0)dF_U(0)e_j \widehat{\chi'_+}(0)+dG_U(0)e_j \\
&= \left(\widehat{\K}(0)dF_U(0)+dG_U(0)\right)e_j\\
&=-c_je_j,
\end{align*}
and we exploit the fact that all $c_j\neq0$, and that the $e_j$ form a basis of $\R^n$. Next, we evaluate the scalar product of $\F_{a_j}(0,0)$ with the elements $e_k$ of the kernel of the adjoint $\cl^*$:
\bqs
\int_\R \langle \F_{a_j}(0,0),e_k \rangle dx= \langle \widehat{\K}(0)dF_U(0)+  dG_U(0) e_j,e_k \rangle = -c_j\delta_{j,k}
\eqs
which, for fixed $j$, is nonzero for $j=k$. Hence, the partial derivative $\F_{a_j}(0,0)$ are not in the range of $\cl$. This proves the lemma.
\end{Proof}

We can now solve \eqref{eq:functionF} with the Implicit Function Theorem and obtain unique solutions $(a,W)(b;\epsilon)$ and thus a solution $U$ of the form \eqref{ansatzappli} to  \eqref{eq:nvclstat}. As outlined in \cite{sandstedescheel:08}, the physically interesting quantity is the jump $U(\infty)-U(-\infty)=a(b;\epsilon)-b$. A straightforward expansion in $\epsilon$ gives
\bqs
U(\infty)-U(-\infty)=a(b;\epsilon)-b=-\epsilon \int_\R \left(\widehat{\K}(0)dF_U(0)+ dG_U(0)\right)^{-1}H(x,0,0)dx+\mathcal{O}(\epsilon^2)
\eqs
which is independent of $b$ to leading order.

The preceding analysis also allows us  to study the case where precisely one characteristic speed $c_{j_0}$ vanishes. In this situation we may further assume that $\langle \widehat{\K}'(0)dF_U(0) e_{j_0}, e_{j_0} \rangle\neq 0$, such that $\nu=0$ is a simple zero of $\det(\widehat{\K}(\nu)dF_U(0)+dG_U(0))=0$ and $(\widehat{\K}(0)dF_U(0)+dG_U(0))e_{j_0}=0$. We directly see that the Fredholm index of $\widetilde{\cl}$ and thus $\cl$ in $L^2_\eta(\R,\R^n)$ is now $-(n+1)$, since $\nu=0$ has multiplicity $n+1$ as a solution of \eqref{eq:charalin}. The kernel of the adjoint operator $\cl^*$ is spanned by the constant functions $e_j$ and the linear function $x e_{j_0}$. Indeed, we have
\begin{align*}
\cl^*(xe_{j_0})&=-dF_U(0)^t \K_-^t\ast e_{j_0}-dG_U(0)e_{j_0}\\
&=-\left(dF_U(0)^t \widehat{\K}_-^t(0)+dG_U(0)^t \right)e_{j_0}\\
&=-\left(\widehat{\K}(0)dF_U(0)+dG_U(0) \right)e_{j_0}\\
&=0.
\end{align*}

We can once again use the ansatz \eqref{ansatzappli} and arrive at the function $\F$ given in \eqref{eq:functionF}.

\begin{lem}
Assume that $\widehat{\K}(0)dF_U(0)+dG_U(0)$ has distinct real eigenvalues with a simple eigenvalue at $\nu=0$ with eigenvector $e_{j_0}$. We also suppose that $\langle \widehat{\K}'(0)dF_U(0) e_{j_0}, e_{j_0} \rangle\neq 0$. Then the linearization of $\F$ with respect to $(a,b_{j_0},W)$ is invertible at $(0;0)$.
\end{lem}

\begin{Proof}
One readily verifies that the partial derivatives with respect to $(a_j)_{j=1,\dots,n}$ and $b_{j_0}$ are linearly independent and that for each fixed $j=1,\dots,n$, $j\neq j_0$, we have
\bqs
\int_\R \langle \F_{a_j}(0,0),e_k \rangle dx= \langle \widehat{\K}(0)dF_U(0)+  dG_U(0) e_j,e_k \rangle = -c_j\delta_{j,k},
\eqs
which is non zero for $j=k$. Lastly, 
\begin{align*}
\int_\R \langle \F_{a_{j_0}}(0,0),x e_{j_0} \rangle dx&=-\int_\R \langle \left(\K \ast (dF_U(0)\chi_+)+dG_U(0)\chi_+ \right)e_{j_0},e_{j_0}\rangle dx\\
&=-\langle \widehat{\K}'(0)dF_U(0) e_{j_0}, e_{j_0} \rangle\neq 0
\end{align*}
and similarly
\bqs
\int_\R \langle \F_{b_{j_0}}(0,0),x e_{j_0} \rangle dx=\langle \widehat{\K}'(0)dF_U(0) e_{j_0}, e_{j_0} \rangle\neq 0
\eqs
so that $\F_{a_{j_0}}(0;0)$ and $\F_{b_{j_0}}(0;0)$ do not lie in the range of $\F_W(0;0)$. Thus $\F_{a,b_{j_0},W}(0;0)$ is invertible.
\end{Proof}

We can therefore solve \eqref{eq:functionF} using the Implicit Function Theorem and obtain a unique solution $(a,b_{j_0},W)$ as functions of $\left((b_j)_{j=1,\dots,n,~ j \neq j_0};\epsilon\right)$. In that case we have that the solution $U$ selects both $a_{j_0}$ and $b_{j_0}$ via
\bqs
a_{j_0}=M\epsilon+\mathcal{O}(\epsilon^2),\quad b_{j_0}=-M\epsilon+\mathcal{O}(\epsilon^2),\quad M:=\int_\R \frac{x\langle H(x,0,0),e_{j_0}\rangle}{\langle \widehat{\K}'(0)dF_U(0) e_{j_0}, e_{j_0} \rangle}dx.
\eqs
When $M\neq 0$, the difference between the number of positive characteristic speeds at $\infty$ and $-\infty$ is two, and the viscous profile is a Lax shock or under compressive shock of index $2$. 

Summarizing, we have shown that nonlocal conservation laws behave in a very similar fashion as local conservation laws when subject to local source terms. Sources that move with non-characteristic speed cause a jump across the inhomogeneity, while number of ingoing and outgoing characteristics are equal. Sources that move with characteristic speed are able to act as sources with respect to the characteristic speed, so that the number of outgoing characteristics exceeds the number of incoming characteristics by two. 

In both cases, stationary profiles are \emph{smooth}, similar to what one would expect from a viscous conservation law. Loosely speaking, smoothing here is provided by dispersal through the nonlocal term rather than smoothing by viscosity. 

\subsection{Edge bifurcations and the nonlocal Gap Lemma}\label{secGL}

We show how our methods can be used to study eigenvalue problems near the edge of the essential spectrum. Motivated most recently by questions on stability of coherent structures, such as solitons in dispersive equations and viscous shock profiles, there has been significant interest in studying spectra of operators near the edge of the essential spectrum. In the original works \cite{gz,ks}, a Wronskian-type function that tracks eigenvalues and multiplicities via its roots was extended into the essential spectrum, exploiting the fact that coefficients of the linearized problem converge exponentially as $|x|\to\infty$. While Wronskians are usually finite-dimensional, extensions are sometimes possible to infinite-dimensional systems, using exponential dichotomies and Lyapunov-Schmidt reduction to obtain reduced Wronskians. 

Gap Lemma type arguments had been used routinely in the theory of Schr\"odinger operators, providing extensions of scattering coefficients into and across the continuous spectrum. One is often interested in tracking how eigenvalues may emerge out of the essential spectrum when parameters are varied. It was observed early that small localized traps inserted into a free Schr\"odinger equation will create bound states in dimensions $n\leq 2$; see \cite{simon:76}. The bound state corresponds to an eigenvalue emerging from the edge of the continuous spectrum. 

We show here how a result analogous to \cite{simon:76} can be proved for nonlocal eigenvalue problems. 
We therefore consider the system
\bqq
\label{eq:GL}
\mathcal{T}(\lambda,\epsilon)\cdot U:=U_\xi+\left(\K+ \epsilon \widetilde{\K}_\xi \right)\ast U -\lambda BU=0, \quad U\in\R^n.
\eqq
Here, $\K,\widetilde{\K}_\xi \in L^1_{\eta_0}(\R,\M_n(\R))$, $B \in \M_n(\R)$, and $\widetilde{\K}_\xi\underset{\xi \rightarrow \pm \infty}{\longrightarrow}0$ in $L^1_{\eta_0}(\R,\M_n(\R))$ such that there exist constants $C>0$ and $\delta>0$ with
\bqs
\left\| \widetilde{\K}(\zeta;\xi) \right\|_n \leq C e^{-\delta |\xi|}, \quad \forall\, \zeta\in\R.
\eqs
We think of \eqref{eq:GL} as coming from a higher-order differential operator such as $\partial_{\xi\xi}$, including nonlocal terms, after rewriting the eigenvalue problem as a first-order system of (nonlocal) differential equations in $\xi$.

\begin{prop}\label{gaplemma}
We assume that the dispersion relation 
\bqs
d(\nu,\lambda)=\det\left(\nu \mathbb{I}_n+\widehat{\K}(\nu) -\lambda B \right)
\eqs
is \emph{diffusive} near $\lambda=0$:
\begin{enumerate}
\item $d(0,0)=d_\nu(0,0)=0$; 
\item $d_{\nu\nu}(0,0) \cdot d_{\lambda}(0,0)<0$; and
\item $d(i\ell,0)\neq0$ for all $\ell\in\R$, $\ell\neq0$.
\end{enumerate}
We also assume that the localized perturbation is generic:
\bqs
M:=\frac{\left\langle \widetilde{\K}_\xi e_0,e_0^*\right\rangle_{L^2(\R,\R^n)}}{\left\langle 2(\mathbb{I}_n+\partial_\nu\widehat{\K}(0))e_1+\partial_{\nu\nu} \widehat{\K}(0)e_0,e_0^* \right\rangle_{\R^n}}\sqrt{-\frac{d_{\nu\nu}(0,0)}{2d_\lambda(0,0)}}\neq 0.
\eqs
Then there exists $\epsilon_0>0$, such that for all $0<M\epsilon<\epsilon_0$ there exist $0\neq U_\epsilon\in H^1(\R,\R^n)$ and $\lambda_*(\epsilon)>0$ so that
\bqs
\mathcal{T}(\lambda_*(\epsilon),\epsilon)\cdot U_\epsilon=0.
\eqs
We also have the asymptotic expansion:
\bqq
\label{eq:asymptexp}
\underset{\epsilon \rightarrow 0^+}{\lim}\frac{\lambda_*(\epsilon)}{\epsilon^2}=M^2.
\eqq
\end{prop}

We prepare the proof of this proposition by reformulating the eigenvalue problem as a nonlinear equation that can be solved with the Implicit Function Theorem near a trivial solution. We first introduce $\lambda=\gamma^2$, so that the dispersion relation has local analytic roots $\gamma \longmapsto \nu_\pm(\gamma) \in \C$. Expanding $d(\nu,\gamma^2)$ in $\gamma^2$, we arrive at the expansion
\bqs
d(\nu,\gamma^2)=\nu^2\frac{d_{\nu\nu}(0,0)}{2}+\gamma^2d_\lambda(0,0)+\mathcal{O}\left(|\nu|^3+|\gamma|^3 \right),
\eqs 
so that to leading order we have
\bqs
\nu_\pm(\gamma)=\pm \sqrt{-\frac{2d_\lambda(0,0)}{d_{\nu\nu}(0,0)}}\gamma+\mathcal{O}(\gamma^2).
\eqs
Associated with these roots can be analytic vectors in the kernel, $\gamma \longmapsto e_\pm(\gamma) \in \C^n$, with
\bqq
\label{analfunc}
\left(\nu_\pm(\gamma)\mathbb{I}_n+\widehat{\K}(\nu_\pm(\gamma)) -\gamma^2B\right)e_\pm(\gamma)=0,
\eqq
and $e_0=e_\pm(0)\neq 0$ solves  $\widehat{\K}(0)e_0=0$.

Following the analysis of the previous section, there exists $\eta_*>0$ such that for each fixed $\eta$ with $0<\eta<\eta_*$, the linear operator $\cl$
\bqs
\cl :U \longmapsto \frac{d}{d\xi}U+\K \ast U,
\eqs
defined on $L^2_\eta(\R,\R^n)$, is Fredholm with index $-2$ and has trivial null space. Indeed, from the above properties, we see that
\bqs
d(\nu,0)=\det\left(\nu \mathbb{I}_n+\widehat{\K}(\nu) \right)=\nu^2\widetilde{d}(\nu), \quad \widetilde{d}(0)\neq 0,
\eqs
with $d(i\ell,0)\neq0$ for all $\ell\in\R$, $\ell\neq0$. This implies that $\nu=0$ is a root with multiplicity $2$ and all other roots have nonzero real part. Thus the Fredholm index of $\cl$ is $-2$ and it is straightforward to check that the kernel of $\cl$ in the exponentially weighted space $L^2_\eta(\R,\R^n)$ is trivial. Thus the kernel of the $L^2$-adjoint $\cl^*$ of $\cl$ considered on $L^2_{-\eta}(\R,\R^n)$ is two-dimensional. Here, the adjoint $\cl^*$ is given via
\bqs
\cl^*:U\longmapsto -\frac{d}{d\xi}U+\K^t_-\ast U,
\eqs
where $\K_-^t(\xi)=\K^t(-\xi)$ for all $\xi\in\R$. Note that
\bqs
\det\left(\widehat{\cl^*}(\nu)\right)=\det\left(-\nu \mathbb{I}_n+\widehat{\K}^t(-\nu) \right)=d(-\nu,0)=\nu^2\widetilde{d}(-\nu),
\eqs
so that there exists $e_0^*\in\R^n$ with $\widehat{\K}^t(0)e_0^*=0$ and thus $\cl^*(e_0^*)=0$. As $d_\nu(0,0)=0$, the following scalar product vanishes:
\bqq
\label{compcond1}
\left\langle (\mathbb{I}_n+\partial_\nu\widehat{\K}(0))e_0,e_0^* \right\rangle_{\R^n} = 0,
\eqq
which ensures the existence of $e_1^*\in\R^n$ so that
\bqq
\label{compcond2}
-\left(\mathbb{I}_n+\partial_\nu\widehat{\K}^t(0)\right)e_0^*+\widehat{\K}^t(0)e_1^*=0.
\eqq
Indeed, the above equation can be solved if $\left\langle\left(\mathbb{I}_n+\partial_\nu\widehat{\K}^t(0)\right)e_0^*,e_0\right\rangle_{\R^n} = 0$, which holds true because of \eqref{compcond1}. We now claim that $\xi e_0^*+e_1^*$ belongs to the kernel of $\cl^*$:
\begin{align*}
\cl^*\left(\xi e_0^*+e_1^*\right)&=\left[-e_0^*+\K^t_-\ast\left(\xi e_0^* \right)\right] +\widehat{\K}^t(0)e_1^*\\
&=\left[-e_0^*-\partial_\nu \widehat{\K}^t(0) e_0^*\right] +\widehat{\K}^t(0)e_1^*\\
&=0.
\end{align*}
Summarizing, the kernel of $\cl^*$, considered on $L^2_{-\eta}(\R,\R^n)$, is spanned by the functions  $e_0^*$ and $\xi e_0^*+e_1^*$.

In the same way, we also define $e_1\in\R^n$ via
\bqq
\label{compcond3}
\left(\mathbb{I}_n+\partial_\nu \widehat{\K}(0) \right)e_0+\widehat{\K}(0)e_1=0.
\eqq
Furthermore, differentiating \eqref{analfunc} with respect to $\gamma$ and evaluating at $\gamma=0$ we obtain
\bqs
\pm \sqrt{-\frac{2d_\lambda(0,0)}{d_{\nu\nu}(0,0)}} \left(\mathbb{I}_n+\partial_\nu \widehat{\K}(0) \right)e_0+\widehat{\K}(0)e_\pm'(0)=0.
\eqs
We see from the above equation and \eqref{compcond3} that $e_\pm'(0)=\pm \sqrt{-\frac{2d_\lambda(0,0)}{d_{\nu\nu}(0,0)}}e_1$. Moreover, combining equations \eqref{compcond2} and \eqref{compcond3} we have the equality
\bqq
\label{compcond4}
\left\langle (\mathbb{I}_n+\partial_\nu\widehat{\K}(0))e_1,e_0^* \right\rangle_{\R^n} = -\left\langle (\mathbb{I}_n+\partial_\nu\widehat{\K}(0))e_0,e_1^* \right\rangle_{\R^n}.
\eqq

The fact that $d_{\nu\nu}(0,0) \neq 0$ ensures that the following quantity is not vanishing:
\bqq
\label{compcond5}
\left\langle (\mathbb{I}_n+\partial_\nu\widehat{\K}(0))e_1,e_0^* \right\rangle_{\R^n}+ \frac{1}{2} \left\langle \partial_{\nu\nu} \widehat{\K}(0)e_0,e_0^*\right\rangle_{\R^n}\neq 0.
\eqq

To find solutions of the eigenvalue problem \eqref{eq:GL}, for small $\epsilon$, we make the following ansatz
\bqq
\label{ansatzGL}
U(\xi)=a_+e_+(\gamma)\chi_+(\xi)e^{\nu_+(\gamma)\xi}+a_-e_-(\gamma)\chi_-(\xi)e^{\nu_-(\gamma)\xi}+w(\xi),
\eqq
where $a_+,a_-\in\R$ and $w\in L^2_\eta(\R,\R^n)$. Here $\chi_+(\xi)=\frac{1+\rho(\xi)}{2}$, where $\rho\in\mathcal{C}^\infty(\R)$ is a smooth even function satisfying $\rho(\xi)=-1$ for all $\xi\leq -1$, $\rho(\xi)=1$ for all $\xi\geq 1$ and $\chi_-(\xi)=1-\chi_+(\xi)$. Substituting the ansatz into \eqref{eq:GL}, we obtain an equation of the form
\bqq
\label{FGL}
\F(a,\gamma,w;\epsilon)=0, \quad \F(~\cdot~;\epsilon):\R^2\times\R\times\R^n\times \mathcal{D}(\cl) \longrightarrow L^2_\eta(\R,\R^n)
\eqq
for $\mathbf{a}=(a_+,a_-)$. We have that $\F((1,1),0,0;0)=0$. For small enough $\eta$, following the analysis conducted in \cite{pogan-scheel:10} and exploiting the localization of $\widetilde{\K}_\xi$, we have that $\F$ is a smooth map. Its linearization at $(\mathbf{a},\gamma,w)=(\mathbf{1},0,0)$ (here for convenience we have denoted $\mathbf{1}=(1,1)$) is given by
\begin{align*}
\F_w(\mathbf{1},0,0;0)&=\cl, \\
\F_{a_\pm}(\mathbf{1},0,0;0)&=\cl\left( \chi_\pm e_0\right), \\
\F_{\gamma}(\mathbf{1},0,0;0)&= \sqrt{-\frac{2d_\lambda(0,0)}{d_{\nu\nu}(0,0)}}\left[\cl\left(\chi_+e_1\right)+\cl\left(\xi \chi_+e_0 \right)\right]-\sqrt{-\frac{2d_\lambda(0,0)}{d_{\nu\nu}(0,0)}}\left[\cl\left(\chi_-e_1\right)+\cl\left(\xi \chi_-e_0 \right)\right]
\end{align*}
where $\F_\mathbf{a}(\mathbf{1},0,0;0)$ and $\F_\gamma(\mathbf{1},0,0;0)$ lie in $L^2_\eta(\R,\R^n)$. 

\begin{lem}\label{l:inve}
Under the above assumptions, the operator
\bqs
\begin{matrix}
\F_{a_-,\gamma,w}(\mathbf{1},0,0;0): &\R\times\R\times L^2_\eta(\R,\R^n)& \longrightarrow & L^2_\eta(\R,\R^n)\\  &(a_-,\gamma,w)& \longmapsto &\F_{a_-}(\mathbf{1},0,0;0)a_-+\F_{\gamma}(\mathbf{1},0,0;0)\gamma+\F_w(\mathbf{1},0,0;0)w
\end{matrix}
\eqs
is invertible.
\end{lem}

\begin{Proof}
We first recall that the cokernel of $\F_w(0;0)$ is spanned by $e_0^*$ and $\xi e_0^*+e_1^*$. We next evaluate the  functional
\bqs
\cl_0u=\left\langle \cl(ue_0),e_0^*\right\rangle_{L^2(\R,\R^n)},
\eqs
with associated symbol $\widehat{\cl_0}(\nu)=\left\langle \left(\nu \mathbb{I}_n+\widehat{\K}(\nu) \right)e_0,e_0^*\right\rangle_{\R^n}$. We have that $\widehat{\cl_0}(0)=\partial_\nu\widehat{\cl_0}(0)=0$, so that there exists $\mathcal{H} \in L^1_{\eta_0}(\R,\M_n(\R))$ such that $\widehat{\cl_0}(\nu)=\nu^2 \left\langle \widehat{\mathcal{H}}(\nu)e_0,e_0^*\right\rangle_{\R^n}= \left\langle  \widehat{\frac{d^2}{d\xi^2}\mathcal{H}}(\nu)e_0,e_0^*\right\rangle_{\R^n}$ with $2\widehat{\mathcal{H}}(0)=\partial_{\nu\nu}\widehat{\K}(0)$. We can rewrite $\cl_0u$ as
\bqs
\cl_0u=\left\langle  \mathcal{H} \ast \left(\frac{d^2}{d\xi^2} u e_0\right),e_0^* \right\rangle_{L^2(\R,\R^n)}.
\eqs
It is now a straightforward computation to evaluate the following quantities:
\begin{align*}
 \cl_0\chi_- &= \left\langle \frac{d^2}{d\xi^2} \mathcal{H} \ast \left( \chi_- e_0\right),e_0^* \right\rangle_{L^2(\R,\R^n)} = 0,\\
\cl_0(\xi \chi_\pm)&=\left\langle \frac{d^2}{d\xi^2} \mathcal{H}\ast \left( \xi \chi_\pm e_0\right),e_0^*\right\rangle _{L^2(\R,\R^n)}\\
& = \left\langle \widehat{\mathcal{H}}(0)e_0,e_0^*\right\rangle_{\R^n} \left(\underset{\xi\rightarrow +\infty }{\lim}\left[\frac{d}{d\xi}\left(\xi\chi_\pm(\xi)\right)\right] -\underset{\xi\rightarrow -\infty }{\lim}\left[\frac{d}{d\xi}\left(\xi\chi_\pm(\xi)\right)\right]  \right)\\
&=\pm\frac{1}{2}\left\langle \partial_{\nu\nu} \widehat{\K}(0)e_0,e_0^*\right\rangle_{\R^n}.
\end{align*}

We can also define the functional
\bqs
\cl_1u=\left\langle \cl(u~e_1),e_0^*\right\rangle_{L^2(\R,\R^n)}
\eqs
such that $\widehat{\cl_1}(\nu)=\left\langle \left(\nu \mathbb{I}_n+\widehat{\K}(\nu) \right)e_1,e_0^*\right\rangle_{\R^n}$ and $\widehat{\cl_1}(0)=0$. Thus, we can find $\mathcal{H}_1 \in L^1_{\eta_0}(\R,\M_n(\R))$ such that $\widehat{\cl_1}(\nu)=\nu \left\langle \widehat{\mathcal{H}}_1(\nu)e_1,e_0^*\right\rangle_{\R^n}= \left\langle  \widehat{\frac{d}{d\xi}\mathcal{H}}_1(\nu)e_1,e_0^*\right\rangle_{\R^n}$ with $\widehat{\mathcal{H}}_1(0)=\mathbb{I}_n+\partial_{\nu}\widehat{\K}(0)$. Using \eqref{compcond4} we find that
\begin{align*}
\cl_1\chi_\pm&= \left\langle \frac{d}{d\xi} \mathcal{H}_1 \ast \left( \chi_- e_1\right),e_0^* \right\rangle_{L^2(\R,\R^n)}\\
&=\pm\left\langle (\mathbb{I}_n+\partial_\nu\widehat{\K}(0))e_1,e_0^* \right\rangle_{\R^n}\\
&=\mp \left\langle (\mathbb{I}_n+\partial_\nu\widehat{\K}(0))e_0,e_1^* \right\rangle_{\R^n}.
\end{align*}

We have thus shown that
\begin{align*}
\left\langle \cl\left(\chi_\pm e_1\right)+\cl\left(\xi \chi_\pm e_0 \right) ,e_0^*\right\rangle_{L^2(\R,\R^n)}&=\cl_1\chi_\pm+ \cl_0(\xi \chi_\pm)\\
&= \mp \left\langle (\mathbb{I}_n+\partial_\nu\widehat{\K}(0))e_0,e_1^* \right\rangle_{\R^n} \pm \frac{1}{2} \left\langle \partial_{\nu\nu} \widehat{\K}(0)e_0,e_0^*\right\rangle_{\R^n}\\
&\neq0.
\end{align*}

Based on similar calculations, we obtain
\begin{align*}
\left\langle\cl(\chi_-e_0),e_1^*+\xi e_0^*\right\rangle_{L^2(\R,\R^n)}&= - \left\langle (\mathbb{I}_n+\partial_\nu\widehat{\K}(0))e_0,e_1^* \right\rangle_{\R^n}+ \frac{1}{2} \left\langle \partial_{\nu\nu} \widehat{\K}(0)e_0,e_0^*\right\rangle_{\R^n}\\
&\neq0.
\end{align*}

Summarizing our results, we have proved that:
\begin{align*}
\left\langle\F_{a_-}(\mathbf{1},0,0;0),e_0^*\right\rangle_{L^2(\R,\R^n)}&=0,\\
\left\langle\F_{a_-}(\mathbf{1},0,0;0),e_1^*+\xi e_0^*\right\rangle_{L^2(\R,\R^n)}&\neq0,\\
\left\langle\F_{\gamma}(\mathbf{1},0,0;0),e_0^*\right\rangle_{L^2(\R,\R^n)}&\neq0.
\end{align*}

Thus $\mathcal{F}_{a_-,\gamma}(0;0)$ span the cokernel of $\cl$, which implies that $\F_{a_-,\gamma,w}(\mathbf{1},0,0;0)$ is invertible, as a Fredholm index 0 operator that is onto.
\end{Proof}
\begin{Proof}[of Proposition \ref{gaplemma}]
Using Lemma \ref{l:inve}, we can solve using the Implicit Function Theorem and obtain a unique solution $(a_-,\gamma,w)$ as a function of $(a_+,\epsilon)$. First, the asymptotic expansion \eqref{eq:asymptexp} follows directly by noticing that, to leading order in $\epsilon$, we have
\bqs
\gamma \left\langle\F_{\gamma}(\mathbf{1},0,0;0),e_0^*\right\rangle_{L^2(\R,\R^n)} +\epsilon \left\langle \widetilde{\K}_\xi e_0,e_0^*\right\rangle_{L^2(\R,\R^n)} + \mathcal{O}(\epsilon^2)=0.
\eqs
Here, we have used the fact that $\left\langle\F_{a_-}(\mathbf{1},0,0;0),e_0^*\right\rangle_{L^2(\R,\R^n)}=\left\langle \cl e_0,e_0^*\right\rangle_{L^2(\R,\R^n)}=0$. Our above computations lead to
\bqs
\left\langle\F_{\gamma}(\mathbf{1},0,0;0),e_0^*\right\rangle_{L^2(\R,\R^n)}=2\sqrt{-\frac{2d_\lambda(0,0)}{d_{\nu\nu}(0,0)}}\left\langle (\mathbb{I}_n+\partial_\nu\widehat{\K}(0))e_1+\frac{1}{2}\partial_{\nu\nu} \widehat{\K}(0)e_0,e_0^* \right\rangle_{\R^n} \neq0.
\eqs
This gives the desired expansion  \eqref{eq:asymptexp} and implies that $\gamma=-M\epsilon +\mathcal{O}(\epsilon^2)$ is of negative sign for $M\epsilon>0$. In order to find have an eigenvalue $\lambda_*(\epsilon)>0$ for \eqref{eq:GL}, we need to check that $U_\epsilon(\xi)$ given in the ansatz \eqref{ansatzGL} belongs to $L^2(\R,\R^n)$. For small $M\epsilon>0$, we have that $\nu_\pm(\gamma)=\mp\sqrt{-\frac{2d_\lambda(0,0)}{d_{\nu\nu}(0,0)}}M\epsilon+\mathcal{O}(\epsilon^2)$, such that $\mp\Re\left(\nu_\pm(\gamma)\right)>0$ and $U_\epsilon$ is exponentially localized. Since for $\lambda>0$, there are no roots $\nu\in i \R$, we know that $\T(\lambda,\epsilon)$ is Fredholm with index zero. Together, this implies that $\T(\lambda,\epsilon)$ possesses a kernel for $\lambda=\lambda_*(\epsilon)$. This completes the proof of Proposition \ref{gaplemma}.

\end{Proof}

\begin{rmk}
Following  \cite[Prop. 5.11]{pogan-scheel:10}, one can show uniqueness and simplicity of the eigenvalue $\lambda_*(\epsilon)$ for $M \epsilon>0$. Also, the analysis here gives a natural extension of the eigenvalue concept into the essential spectrum: for $M \epsilon <0$, we can track the eigenvalue $\lambda_*(\epsilon)$  in smooth fashion as a resonance pole, that is, a function with particular prescribed exponential growth. In this sense, our method here provides an alternative to the Gap Lemma \cite{gz,ks}, were this  possibility of tracking eigenvalues into the essential spectrum was the main objective. 
\end{rmk}

{\noindent \bf Acknowledgments:} GF was partially supported by the National Science Foundation through grant NSF-DMS-1311414. AS was partially supported by the National Science Foundation through grant NSF-DMS-0806614.

\bibliography{plain}

\end{document}